\documentclass{gtart_h}
\def\ifplaintex{\expandafter\ifx\csname documentclass\endcsname\relax}

\def\gtp{{\mathsurround=0pt\it $\cal G\mskip-2mu$eometry \&\ 
$\cal T\!\!$opology $\cal P\!$ublications}}  % GT publications

\def\recd{{\small Received:\qua\receiveddate\ifx\reviseddate\relax
\else\qquad Revised:\qua\reviseddate\fi\par}} 

\def\lognumber#1{\def\thelognumber{#1}}
\def\volumenumber#1{\def\thevolumenumber{#1}}
\def\volumeyear#1{\def\thevolumeyear{#1}}
\def\papernumber#1{\def\thepapernumber{#1}}
\def\pagenumbers#1#2{\def\startpage{#1}\def\finishpage{#2}}
\def\published#1{\def\publishdate{#1}}

\def\received#1{\def\receiveddate{#1}}
\def\revised#1{\def\reviseddate{#1}}
\def\accepted#1{\def\accepteddate{#1}}
\def\asciititle#1{\def\theasciititle{#1}}

\def\asciiaddress#1{\def\theasciiaddress{#1}}

\def\asciikeywords#1{\def\theasciikeywords{#1}}

\let\\\par\let\thelognumber\relax\let\thevolumenumber\relax
\let\thepapernumber\relax\let\thevolumeyear\relax\let\startpage\relax
\let\finishpage\relax\let\publishdate\relax\let\receiveddate\relax
\let\reviseddate\relax\let\accepteddate\relax\let\theasciititle\relax
\let\theasciiauthors\relax\let\theasciiaddress\relax
\let\theasciiabstract\relax\let\theasciikeywords\relax

\let\theasciiemail\relax

\ifplaintex
\font\logobig=cmssbx10 scaled 3836
\font\logomed=cmssbx10 scaled 2557
\else
\font\logobig=cmssbx10 scaled 4200
\font\logomed=cmssbx10 scaled 2800
\fi

\long\def\makeagttitle{   %%% start of definition of \makeagttitle
\count0=\startpage
\agt\hfill      %   Journal title (top left) 
\hbox to 45truept{\vbox to 0pt{\vglue -13truept{\logomed A\kern -.37em{\logobig 
T}\kern -.38em G}\vss}\hss}
\break
{\small Volume \thevolumenumber\ (\thevolumeyear)
\startpage--\finishpage\nl
Published: \publishdate}

\vglue .25truein

{\parskip=0pt\leftskip 0pt plus
1fil\def\\{\par\smallskip}{\Large\bf\thetitle}\par\medskip} \vglue
0.05truein

{\parskip=0pt\leftskip 0pt plus 1fil\def\\{\par}{\sc\theauthors}
\par\medskip}%
 
\vglue 0.03truein 

{\small\leftskip 25truept\rightskip 25truept{\bf Abstract}\stdspace\theabstract

{\bf AMS Classification}\stdspace\theprimaryclass
\ifx\thesecondaryclass\relax\else; \thesecondaryclass\fi\par
{\bf Keywords}\stdspace \thekeywords\par}\vglue 7truept

}   %%%% end of definition of \makeagttitle

\ifplaintex
\font\phead=cmsl9 scaled 950
\font\pnum=cmbx10 scaled 913
\font\pfoot=cmsl9 scaled 950
\headline{\vbox to 0pt{\vskip -4.5mm\line{\small\phead\ifnum
\count0=\startpage ISSN 1472-2739 (on-line) 1472-2747 (printed)
\hfill {\pnum\folio}\else\ifodd\count0\def\\{ }% 
\ifx\theshorttitle\relax\thetitle\else\theshorttitle\fi\hfill{\pnum\folio}
\else\def\\{ and }{\pnum\folio}\hfill\ifx\theshortauthors\relax\theauthors
\else\theshortauthors\fi\fi\fi}\vss}}
\footline{\vbox to 0pt{\vglue 0mm\line{\small\pfoot\ifnum\count0=\startpage
\copyright\ \gtp\hfill\else
\agt, Volume \thevolumenumber\ (\thevolumeyear)\hfill\fi}\vss}}
\else
\font=cmsl9 scaled 1050
\font\lnum=cmbx10 
\font=cmsl9 scaled 1050
\makeatletter
\def\@oddhead{{\small\lhead\ifnum\count0=\startpage ISSN 1472-2739 
(on-line) 1472-2747 (printed)\hfill {\lnum\number\count0}\else\ifodd\count0
\def\\{ }\ifx\theshorttitle\relax \thetitle \else\theshorttitle\fi\hfill
{\lnum\number\count0}\else\def\\{ and }{\lnum\number\count0}
\hfill\ifx\theshortauthors\relax 
\theauthors\else\theshortauthors\fi\fi\fi}}\def\@evenhead{\@oddhead}
\def\@oddfoot{\small\lfoot\ifnum\count0=\startpage\copyright\ \gtp\hfill\else
\agt, Volume \thevolumenumber\ (\thevolumeyear)\hfill\fi}
\def\@evenfoot{\@oddfoot}
\makeatother
\fi
\let\maketitlepage\makeagttitle
\let\makeshorttitle\maketitlepage
\let\maketitle\maketitlepage

\newwrite\gtoutfile
\long\gdef\makeheadfile{  %%% start of definition of \makeheadfile
{\def\\{, }\def\s{ }
\immediate\openout\gtoutfile head.xxx
\immediate\write\gtoutfile{Proxy-for: \ifx\theasciiauthors\relax
\theauthors\else\theasciiauthors\fi\s<\ifx\theasciiemail\relax\theemail\else\theasciiemail\fi>}
\immediate\write\gtoutfile{\noexpand\\}
\immediate\write\gtoutfile{Authors: \ifx\theasciiauthors\relax
\theauthors\else\theasciiauthors\fi}
{\def\\{ }\immediate\write\gtoutfile{Title: \ifx\theasciititle\relax
\thetitle\else\theasciititle\fi}}
\immediate\write\gtoutfile{Subj-class: GT or SG, GR etc}
\immediate\write\gtoutfile{MSC-class: \theprimaryclass\ifx\thesecondaryclass\relax\else, \thesecondaryclass\fi}
\immediate\write\gtoutfile{Journal-ref: Algebr. Geom. Topol. \thevolumenumber\s
(\thevolumeyear) \startpage-\finishpage}
\immediate\write\gtoutfile{Comments: Published by Algebraic and
Geometric Topology at}
\immediate\write\gtoutfile{\s\s\s  http://www.maths.warwick.ac.uk/agt/AGTVol\thevolumenumber/agt-\thevolumenumber-\thepapernumber.abs.html}
\immediate\write\gtoutfile{\noexpand\\}
\immediate\write\gtoutfile{}
\ifx\theasciiabstract\relax
\immediate\write\gtoutfile{\theabstract}\else
\immediate\write\gtoutfile{\theasciiabstract}\fi
\immediate\write\gtoutfile{}
\immediate\write\gtoutfile{\noexpand\\}
\immediate\write\gtoutfile{}
\immediate\closeout\gtoutfile}}  %%% end of definition of \makeheadfile

\def\maketitlepage{\makeagttitle\makeheadfile}
\let\makeshorttitle\maketitlepage
\let\maketitle\maketitlepage

\lognumber{34}
\volumenumber{5}
\volumeyear{2005}
\papernumber{34}
\pagenumbers{785}{833}
\received{27 October 2003} 
\revised{14 July 2005}
\accepted{24 January 2005}
\published{24 July 2005}

\usepackage{amsmath,amssymb,xspace,psfrag,graphicx}

\newcommand{\R}{\mathbb R}
\newcommand{\N}{\mathbb N}
\newcommand{\Q}{\mathbb Q}
\newcommand{\Z}{\mathbb Z} 

\newtheorem{theorem}{Theorem}[section]
\newtheorem{lemma}[theorem]{Lemma}
\newtheorem{cor}[theorem]{Corollary}
\newtheorem{prop}[theorem]{Proposition}

\theoremstyle{definition}
\newtheorem{dfn}[theorem]{Definition}

\newtheorem*{rem}{Remark}

\numberwithin{figure}{section}

\begin{document}

\title{Tight contact structures on Seifert manifolds\\over $T^2$ with one 
singular fibre}
\shorttitle{Tight contact structures on Seifert manifolds over $T^2$}
\asciititle{Tight contact structures on Seifert manifolds over T^2 with one 
singular fibre}
 
\author{Paolo Ghiggini} 

\address{Dipartimento di Matematica ``L.~Tonelli'', Universit\`a di 
Pisa\\Largo Pontecorvo 5, I--56127 Pisa, Italy}
\asciiaddress{Dipartimento di Matematica `L. Tonelli', Universita di 
Pisa\\Largo Pontecorvo 5, I--56127 Pisa, Italy}

\email{ghiggini@mail.dm.unipi.it}

\begin{abstract}
In this article we classify up to isotopy tight contact structures on 
Seifert manifolds over the torus with one singular fibre. 
\end{abstract}
\primaryclass{57R17}
\secondaryclass{57M50}
\keywords{Contact structure, tight, Seifert $3$--manifold, convex surface}
\asciikeywords{Contact structure, tight, Seifert 3-manifold, convex surface}

\maketitle

\section{Introduction}  
  A {\em contact structure} on a $3$--manifold $M$ is a
tangent $2$--plane field $\xi$ which is 
the kernel of a differentiable $1$--form $\alpha$ such that $\alpha
\land d \alpha$ is a nowhere vanishing $3$--form. Contact structures 
on $3$--manifolds split into two families.
A contact structure $\xi$ is {\em overtwisted} if there exists an
embedded disc $D \subset M$ such that 
$TD|_{\partial D} \equiv \xi|_{\partial D}$. A contact structure is 
{\em tight} if it is not overtwisted.
The disc $D$ is called, with an abuse of terminology, an overtwisted disc.

Overtwisted contact structures are much more common and flexible objects than 
the tight ones, in fact any $3$--manifold admits an overtwisted contact 
structure and on a closed $3$--manifold two 
overtwisted contact structures are isotopic if and only if they are homotopic 
as plane fields (Eliashberg \cite{eliashberg:4}).
On the contrary, the classification of tight contact structures is still at its
beginning.  
For a survey of contact structures, see
\cite{aebischer,eliashberg:1,etnyre:1,giroux:b}. 

In the last
decade there has been a dramatic growth of the three--dimensional methods in 
contact 
topology starting from the definition of convex surfaces in Giroux's paper
\cite{giroux:1}.  Convex surfaces are the main tool
to perform cut-and-paste operations on contact manifolds. Applying this 
technique, Kanda \cite{kanda:1} and, independently, Giroux, classified the
 tight contact structures on the three--torus $T^3$. Later,
Honda \cite{honda:1} and Giroux \cite{giroux:4} classified
the tight contact structures on lens spaces, the solid torus $D^2 \times S^1$ and 
the thickened torus $T^2 \times I$. In \cite{honda:1}, Honda introduced the notion
 of {\it bypass}, a tool which allows one to handle contact 
topological problems in a combinatorial way (see \cite{honda:1}, Section 
3.4). In this paper we will assume that the reader is familiar with the 
material in \cite{giroux:1} and \cite{honda:1}.
 
The solid torus and the thickened torus can be thought of
 as basic building blocks for a number of other three dimensional 
manifolds. In fact, shortly after, Honda \cite{honda:2} gave a complete 
classification of tight contact structures on  $T^2$--bundles over $S^1$ and 
$S^1$--bundles over surfaces. At the same time Giroux \cite{giroux:5} obtained 
almost complete results on the same manifolds.
  
Other classification results are partial or sporadic. The most important of 
them are  the non existence of tight contact structures on the 
Poincar{\'e} homology sphere with opposite orientation $- \Sigma (2,3,5)$
in \cite{etnyre-honda:1} and the coarse classification which characterises the 
three--manifolds which carry infinitely many tight contact structures, 
\cite{colin:2,CGH:3,CGH:2,honda-kazez-matic:1}. 
A complete classification is also known for the Seifert manifolds over $S^2$
with three singular fibres $\pm \Sigma (2,3, 11)$, 
\cite{ghiggini-schonenberger}. Moreover, there are  partial results on fibred 
hyperbolic
three--manifolds \cite{honda-kazez-matic:2}, which are the only non Seifert
manifolds in the list so far. During the preparation of this article tight
contact structures have been classified on small Seifert manifolds with 
integer Euler class $e_0 \neq -2, -1$, \cite{gls:1,wu:2}.

Our aim  is to give a complete isotopy classification of tight 
contact 
structures on Seifert manifolds over the torus $T^2$ with one singular fibre. 
Fix $e_0 \in \Z$ and $r \in (0,1) \cap \Q$, and let $T(e_0)$ be the circle bundle over
$T^2$ with Euler class $e_0$. We denote by $M(e_0,r)$ the Seifert manifold 
obtained by $(- \frac{1}{r})$--surgery along a fibre of $T(e_0)$. 
 The tight contact structures on  $M(e_0,r)$ and $T(e_0)$
 split into two families, according to 
their behaviour with respect to the finite coverings  induced by a finite 
covering of $T^2$. We will call {\em generic} those tight contact structures 
which remain tight after pulling back to such coverings, and {\em exceptional}
those ones which become overtwisted. 
The set of isotopy classes of generic tight contact structures on $M(e_0,r)$
splits into infinitely many
sub-families parametrised by the isotopy classes of the generic tight contact 
structures on $T(e_0)$. Each sub-family contains finitely many isotopy classes 
of tight contact 
structures which are obtained by Legendrian surgery on the generic tight 
contact structure on $T(e_0)$ labelling the sub-family.

The isotopy classes of exceptional tight contact structures on $M(e_0,r)$
form a finite family, whose cardinality depends on $e_0$ and $r$.
If $e_0 \leq 0$ there are no exceptional tight contact structures on $M(e_0,r)$.
If $e_0 \geq 2$, all exceptional tight contact structures on $M(e_0,r)$ are obtained
by Legendrian surgery on the exceptional tight contact structures on $T(e_0)$
which, however, are not fillable
by \cite{lisca-stipsicz:1}. If $e_0=1$, the exceptional tight contact structures
on $M(e_0,r)$ have
 no tight analogue
on $T(e_0)$.  They are 
obtained by Legendrian surgery
on  overtwisted contact structures and there seems to be no natural way to 
express them
as Legendrian surgery on a tight contact structure. When $e_0=1,2$ the 
exceptional tight contact structures show an unexpected interplay between
the corresponding contact structure on $T(e_0)$ and the 
surgery data. See Theorem \ref{eccezionali}.

\rk{Acknowledgements}
I would like to thank the American Institute of Mathematics and Stanford 
University for their support during the Special Quarter in Contact Geometry
held in Autumn 2000, when this work moved its first steps.
I also thank The University of Georgia at Athens for its support
in the Spring 2002. I am very grateful to Emmanuel Giroux for suggesting  
this problem  and to John Etnyre and especially to Ko Honda for their 
encouragement, many useful discussions, and for helping me to physically 
survive in my first days in Palo Alto during the Contact Geometry Quarter.   
A special thank you to Ko Honda for pointing out some gaps in the 
earlier version of this manuscript. Finally, I thank Riccardo Murri and 
Antonio Messina for their steady computer support.

The author is a member of EDGE, Research Training Network
HPRN--CT--2000--00101, supported by The European Human Potential
Programme.
%%%%%%%%%%%%%%%%%%%%%%%%%%%%%%%%%%%%%%%%%%%%%%%%%%%%%%%%%%%%%%%%%%%%%%%%%%%%%%%
\section{Statement of results}
Let $M$ be an oriented $3$--manifold. The set of isotopy classes of tight 
contact structures on $M$ will be denoted by $\hbox{Tight}(M)$. If $\partial M \neq \emptyset$, 
and ${\mathcal F}$ is a singular foliation on $\partial M$, 
$\hbox{Tight}(M, {\mathcal F})$ will denote the set of tight contact 
structures on $M$ which induce the characteristic foliation ${\mathcal F}$ on 
$\partial M$, modulo isotopies fixed on the boundary. If ${\mathcal F}$ and 
${\mathcal G}$ are two singular foliations on $\partial M$ adapted to the same 
dividing set $\Gamma_{\partial M}$, then $\hbox{Tight}(M, {\mathcal F})$ and 
$\hbox{Tight}(M, {\mathcal G})$ are canonically identified, therefore we will
 write $\hbox{Tight}(M, \Gamma_{\partial M})$ in place of $\hbox{Tight}(M, {\mathcal F})$
for any ${\mathcal F}$ adapted to $\Gamma_{\partial M}$.

Recall that we denote by $T(e_0)$, for $e_0 \in \Z$, the $S^1$--bundle
over $T^2$ with Euler class $e_0$, and by $M(e_0,r)$, for $r \in \Q \cap 
(0,1)$, the Seifert manifold
over $T^2$ obtained by $(- \frac{1}{r})$--surgery along a fibre $R$ of $T(e_0)$.
Here the surgery coefficient is calculated with respect to the standard framing
on $R$. More explicitly, consider a tubular neighbourhood 
$\nu R \subset T(e_0)$ of $R$, and identify $- \partial (T(e_0) \setminus \nu R)$ to $\R^2 / \Z^2$ so that 
$\left ( \begin{array}{c} 1 \\ 0 \end{array} \right )$ is the direction of the 
meridian of $\nu R$ and  $\left ( \begin{array}{c} 0 \\ 1 \end{array} \right )$is 
the direction of the fibres. Then  $M(e_0,r)$ is the manifold obtained by gluing
a solid torus $D^2 \times S^1$ to $T(e_0) \setminus \nu R$ by the map
\[A(r)\co \partial D^2 \times S^1 \to - \partial (T(e_0) \setminus \nu R)\]
represented by the matrix 
\[A(r)= \left ( \begin{matrix}
\alpha & \alpha' \\
- \beta & - \beta' 
\end{matrix} \right ) \in SL(2, \Z)\]
where $r= \frac{\beta}{\alpha}$ and $0 \leq \alpha' < \alpha$.
The image of $\{ 0 \} \times S^1 \subset D^2 \times S^1$ in $M(e_0,r)$ is 
called the {\em singular fibre}. The images of the fibres of $T(e_0)$ are 
called {\em regular fibres}. 
See \cite{fomenko-matveev,hatcher,orlik} for more about 
Seifert manifolds.

Let $M$ be a Seifert manifold, possibly without singular fibres, with non 
simply connected base. Let $R \subset
 M$ be a curve isotopic to a regular fibre. In the following such curve will 
be called a {\em vertical curve}. Following Kanda \cite{kanda:1}, we define the 
{\em canonical framing} of $R$ as the framing induced by any incompressible 
torus $T \subset M$ containing $R$. Unless stated otherwise, the twisting number 
of Legendrian vertical curves will be calculated with respect to the 
canonical framing.
\begin{dfn}
Let $M$ be a Seifert fibred manifold over an oriented non simply 
connected surface.
Given a regular fibre $R \subset M$ and a contact structure $\xi$ on $M$, we define 
the {\em maximal twisting number} of $\xi$ as
\[t(\xi)= \max_{L \in  {\mathcal S}} \min \{ tb(L), 0 \}\]
where ${\mathcal S}$ is the set of all Legendrian curves $L \subset M$ isotopic
to $R$.
\end{dfn}
It is clear that the number $t(\xi)$  does not depend on the choice of $R$, 
and is an
isotopy invariant of $\xi$, therefore it defines a function 
\[t\co \hbox{Tight}(M) \to \Z_{\leq 0}.\]
Seifert fibred manifolds over a surface of genus $g>0$ have a 
distinguished family of coverings: namely, the coverings induced by a covering 
of the base.

\begin{dfn}
A tight contact structure on a Seifert fibred manifold $M$ is of 
{\em generic type} if it remains tight after pull-back with respect to any 
covering of $M$ induced by a finite covering of the base. A tight 
contact structure on $M$ is {\em exceptional} if it becomes overtwisted after 
pull-back  with respect to any covering of $M$ induced by a finite covering of 
the base.
\end{dfn}
\noindent
We denote the set of the isotopy classes of the generic tight contact 
structures on $M$ by $\hbox{\rm Tight}_0(M)$. 

In the following theorem $\Gamma_s$ will be a dividing set 
on $T^2$ with $\# \Gamma_s=2$ and slope $s$.
Every rational number $- \frac{p}{q} <-1 $ has a unique finite continued 
fraction expansion
\[- \frac{p}{q} = d_0 - \frac{1}{d_1 - \frac{1}{\ddots - \frac{1}{d_n}}}\]
with $d_i \leq -2$ for $i>0$. We denote this expansion by $- \frac{p}{q}=[d_0, \ldots ,
d_n]$.

\begin{theorem}\label{numerouno}
All tight contact structures on $M(e_0,r)$ are either of generic type or 
exceptional. There exists a map
\[bg\co \hbox{\rm Tight}_0(M(e_0,r)) \longrightarrow \hbox{\rm Tight}_0(T(e_0))\]
such that, given $\xi_0 \in \hbox{\rm Tight}_0(T(e_0))$,  
\begin{itemize}
\item $bg^{-1}(\xi_0)= \emptyset$ if $t(\xi_0) \leq - \frac{1}{r}$,
\item $bg^{-1}(\xi_0)$ is in natural bijection with
 $\hbox{\rm Tight}(D^2 \times S^1, \, A(r)^{-1} \Gamma_{\frac{1}{t(\xi_0)}})$ and has cardinality
$|(d_0-t(\xi_0))(d_1+1) \ldots (d_n+1)|$, where $[d_0, \ldots ,d_n]$ is the continued fraction 
expansion of $- \frac 1r$, if $t(\xi_0)> - \frac{1}{r}$. 
\end{itemize}

The exceptional tight contact structures exist only when $e_0>0$ and all have 
maximal twisting number $t=0$. Their number is always finite and is
\begin{itemize}
\item $2|(d_0+1) \ldots (d_k+1)|$ if $e_0>2$,
\item $|(d_0-1)(d_1+1) \ldots (d_k+1))|$ if $e_0=2$,
\item $|d_1(d_2+1) \ldots (d_k+1)|$ if $e_0=1$.
\end{itemize}
The last expression has to be interpreted as $2$ when $- \frac{1}{r}=d_0 \in \Z$.
\end{theorem}
The map $bg$ is constructed by removing a tubular neighbourhood of the 
singular fibre $V$ such that $- \partial (M(e_0, r) \setminus V)$ is convex with slope
$\frac{1}{t(\xi)}$ and gluing $D^2 \times S^1$ with the unique tight contact structure 
with boundary slope $\frac{1}{t(\xi)}$ to $- \partial (M(e_0, r) 
\setminus V)$ via the identity map. The identification of $bg^{-1}(\xi_0)$ with
$\hbox{\rm Tight}(D^2 \times S^1, \, A(r)^{-1} \Gamma_{\frac{1}{t(\xi_0)}})$ is 
given by the restriction $(M(e_0,r), \, \xi) \mapsto (V, \, \xi |_V)$.
The fact that the map $bg$ and the restriction $(M(e_0,r), \, \xi) \mapsto (V, \, \xi |_V)$ 
are well defined up to isotopy is part of the statement.  
Theorem \ref{numerouno} exhibits each generic tight contact structure on 
$M(e_0,r)$ as a contact surgery in the sense of \cite{ding-geiges:2} on a 
generic tight contact structure on $T(e_0)$. Moreover, the condition 
$t(\xi_0)> - \frac{1}{r}$ implies that it is a negative contact surgery,
which means that the surgery coefficient, calculated with respect to the 
contact framing, is negative. The expression for the cardinality of $bg^{-1}(\xi_0)$ 
is a consequence of the following lemma, which is simply the 
classification of tight contact structures on solid tori \cite{honda:1}, 
Theorem 2.3 applied to $\xi|_V$ after a change of coordinates. For benefit of the
reader we sketch here how to deduce this lemma from Honda's Theorem.

\begin{lemma}\label{solidi}
$\hbox{\rm Tight}(D^2 \times S^1, \, A(r)^{-1} \Gamma_{\frac{1}{t(\xi_0)}})$ is a nonempty finite 
set with cardinality
\[ |\hbox{\rm Tight}(D^2 \times S^1, \, A(r)^{-1} \Gamma_{\frac{1}{t(\xi_0)}})|=|(d_0-t(\xi_0))(d_1+1) \ldots 
(d_n+1)|, \] 
where $[d_0, \ldots ,d_n]$ is the continued fraction expansion of $- \frac{1}{r}$.
\end{lemma}

\begin{proof}
Let $r'=\frac{1}{\frac{1}{r}+t(\xi)+1}$ so, by a direct check, 
$A(r)^{-1} \Gamma_{\frac{1}{t(\xi_0)}}$ and $A(r')^{-1} \Gamma_{-1}$ have the same 
slope $s'$. By \cite{ding-geiges:2}, proof of Proposition 3, if 
$- \frac{1}{r'}$ has the continued fraction expansion $- \frac{1}{r'}=
[d_0', \ldots ,d_n']$, then $s'$ has the continued fraction expansion 
$s'=[r_n', \ldots ,r_0'+1]$. By \cite{honda:1}, Theorem 2.3, 
\[|\hbox{\rm Tight}(D^2 \times S^1, \, A(r)^{-1} \Gamma_{\frac{1}{t(\xi_0)}})|=|(d_0'+1)(d_1'+1) \ldots 
(d_n'+1)|\] provided that $s'<-1$.
As $d_0'=d_0-(t(\xi)+1)$ and $d_i'=d_i$ for $i>0$, we have $s'<d_n+1 \leq -1$ and 
$|d_0'(d_1'+1) \ldots (d_n'+1)|=|(d_0-t(\xi_0))(d_1+1) \ldots (d_n+1)|$.
\end{proof}
%%%%%%%%%%%%%%%%%%%%%%%%%%%%%%%%%%%%%%%%%%%%%%%%%%%%%%%%%%%%%%%%%%%%%%%%%%%%%%
\section{Tight contact structures on $T(e_0)$}
The tight contact structures on $T(e_0)$ have been classified in \cite{honda:2}
 and in \cite{giroux:4}. The material in this section is taken primarily from 
\cite{honda:2}, adapting statements and notation to our purposes. In order
to fix terminology and notations, we start with a digression about 
characteristic foliations on tori in tight contact manifold before focusing 
on the classification of tight contact structures on $T(e_0)$.

\subsection{Characteristic foliation on tori}
If $T$ is a convex torus in a tight contact manifold $(M, \, \xi)$, by Giroux's 
Tightness Criterion \cite{honda:1} Lemma 4.2, its dividing set $\Gamma_T$ contains no
dividing curve bounding a disc in $T$, therefore it consists of an even number of
 closed, parallel, homotopically non trivial curves. 

\begin{dfn}
If $\gamma$ is a dividing curve of $T$, we call the quantity $s(T)= [\gamma] \in {\mathbb P}
(H_1(T, \Q))$ the slope of the convex torus $T$. 
\end{dfn}

The choice of an 
identification $T \cong \R^2 / \Z^2$ gives an identification ${\mathbb P}(H_1(T, \Q))
\cong \Q \cup \{ \infty \}$, hence we will more often see the slope as a rational number.

\begin{dfn}
We call the {\em division number} of 
$T$ the number $\hbox{div}(T)= \frac{1}{2} \# \Gamma_T$. If $\hbox{div}(T)= 1$ we  say 
that $T$ is {\em minimal}. 
\end{dfn}
 
\begin{figure}[ht!]
\centering
\includegraphics[width=3.5cm]{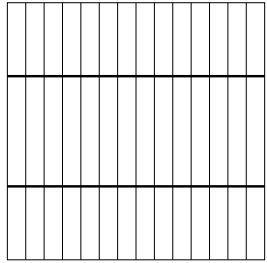}
\caption{Characteristic foliation on a convex torus in standard form with 
vertical Legendrian ruling and two horizontal Legendrian divides}
\label{standard.fig}
\end{figure}

Given a dividing set $\Gamma_T$ on a torus $T$ in a tight contact manifold, there 
is a canonical family of characteristic foliations adapted to $\Gamma_T$. Fix a slope
$r \neq s(T)$ and consider on $T$ the singular foliation consisting of a 
$1$-parameter family of closed curves with slope $r$, called {\em Legendrian 
rulings}, and a closed curve of singularities with slope $s(T)$ called 
{\em Legendrian divide} in each component of $T \setminus \Gamma_T$. See Figure 
\ref{standard.fig} for an illustration. A torus with a characteristic 
foliation of this type is called a convex torus in {\em standard form}, or a
{\em standard torus}.

 As an immediate consequence of Giroux's Flexibility theorem, any convex torus 
$T$ with slope $s(T)$ in a tight contact manifold can be put in standard form
with ruling slope $r$ by a $C^0$-small  perturbation, provided that $r \neq s(T)$.

Sometimes we will need to consider non convex tori of a 
particular kind. 
\begin{dfn}
A {\em pre-Lagrangian torus} is a torus embedded in a contact 
manifold, whose  characteristic foliation after a change of coordinates is 
isotopic to a linear foliation with closed leaves.
\end{dfn}
Suppose we have chosen coordinates on a neighbourhood of a pre-Lagrangian torus
$T$ so that $T= \{ y=0 \}$, and the characteristic foliation of $T$ has slope 
$0$. Then the contact form in a neighbourhood of $T$ is given by $dz-y \, dx$. 
 Pre-Lagrangian tori can be perturbed into convex tori, as explained in the 
following lemma.
\begin{lemma} \label{prelagrangiano}
Let $T$ be a pre-Lagrangian torus whose characteristic foliation has closed 
leaves with slope $s$. Then, for any 
natural number $n>0$, $T$ can be put in standard form with $2n$ dividing
 curves with slope $s$ by a $C^{\infty}$-small perturbation.
\end{lemma}
\begin{proof}
Let $T$ be the given pre-Lagrangian torus. Put coordinates $(x,y,z) \in \R / \Z 
\times I \times \R / \Z$ in a tubular neighbourhood 
$N$ of $T$ such that $T= \{ y=0 \}$ and the contact form is $\alpha =dz-y \, dx$, then
 consider the embedding $i: T^2 \to N$ given by $i: (u,v) \mapsto (u, \epsilon \sin (2 \pi nv),
v)$. After identifying $T^2$ with the image of $i$, the characteristic foliation
is given by the form $i^* \alpha =dv- \epsilon \sin (2 \pi nv)du$. 

Fix the area form $\omega = du \land dv$ on $T^2$, then the characteristic foliation is 
directed by a vector field $X$ such that $\iota_X (\omega)= i^* \alpha$. Since $L_X \omega = d i^* \alpha =
 2 \pi n \epsilon \cos (2 \pi n v) du \land dv$, the set $\Gamma = \{ L_X \omega = 0 \}$ consists of 
$2n$ parallel simple closed curves with slope $0$. The vector field $X$ expands 
$\omega$ where $L_X \omega$ is a positive multiple of $\omega$, and $-X$
expands $\omega$ where $L_X \omega$ is a negative multiple of $\omega$, therefore, by 
\cite{giroux:1} Proposition II.2.1, $\Gamma$ is dividing set for the characteristic 
foliation of $T$.
\end{proof}
 
\subsection{Tight contact structures with $t<0$}
\begin{theorem}[\cite{honda:2}, Lemma 2.7]\label{ciccio}
If $e_0<0$, then on $T(e_0)$ there are $|e_0-1|$ distinct tight
contact structures with $t<0$.
\end{theorem}
By a direct check of the definition of such tight contact structures, see 
\cite{honda:2}, Case 9, it follows that only $2$ of the $|e_0-1|$ are
universally 
tight, but all remain tight if lifted to a covering of $T(e_0)$ induced by a 
finite covering of the base $T^2$. 
\begin{theorem} \label{riempibili}
The tight contact structures with $t<0$ on $T(e_0)$, when $e_0<0$, are Stein
fillable.
\end{theorem}
\begin{proof} In \cite{gompf:1} Gompf constructed $|e_0-1|$ Stein fillings of 
$T(e_0)$ when $e_0 \leq 0$: see \cite{gompf:1},
Figure 36 (c)  for a surgery presentation of the Stein filling of
$T(0)=T^3$. When $e_0 \leq -1$, the Stein fillings of $T(e_0)$ are obtained by 
Legendrian surgery on a stabilisation of the knot in \cite{gompf:1}, Figure 
36(c). 
All the Stein fillings 
obtained in such way are diffeomorphic to the disc bundle over $T^2$ with
Euler class $e_0$, but their complex structures have
different first Chern classes determined by the rotation number of the Legendrian
knot, as shown 
in \cite{gompf:1}, Proposition 2.3. The tight contact structures induced on the
boundary are pairwise non isotopic by \cite{lisca-matic:2}, Corollary 4.2.

To prove Theorem \ref{riempibili}, we need to show  that the $|e_0-1|$ tight 
contact 
structures on $T(e_0)$ induced by the different Stein structures described
above  have $t<0$. 
Let $W$ be the disk bundle over $T^2$ with Euler class $e_0$, and $D \subset W$ a
fibre with Legendrian boundary $\partial D = K$. The 
{\em slice  Thurston--Bennequin invariant} $tb_D(K) \in \Z$ is defined in 
\cite{lisca-matic:2}, Definition 3.1 as the obstruction
to extending the positively oriented normal of the contact structure restricted 
to $K$ to a nowhere vanishing section of the normal bundle of $D$. It can be 
defined equivalently as the twisting number of $K$ computed with respect to the 
framing induced on $K$ by the restriction of a nowhere vanishing section of the 
normal bundle of $D$.
The framing on the normal bundle of $D$ induced by the disc bundle structure
over $W$ restricts to the framing on the normal bundle of $\partial D=K \subset T(e_0)$ 
induced by the circle bundle structure on $T(e_0)$.  

The bundle framing of $K$ coincides with the canonical framing, therefore the 
Thurston--Bennequin number $\hbox{tb}_D(K)$  and the twisting number 
$\hbox{tb}(K)$ defined by the canonical framing coincide. By the slice 
Thurston--Bennequin inequality \cite{lisca-matic:2}, Theorem 3.4 $tb_D(K) \leq -1$ 
for any Legendrian knot $K$ in $(T(e_0), \, \xi)$ smoothly isotopic to a fibre of 
$T(e_0)$, therefore $t(\xi)<0$.
On the other hand there are exactly $|e_0-1|$ tight contact structures on $T(e_0)$
 with $t<0$, so any tight contact structure on $T(e_0)$ with $t<0$ 
must be Stein fillable for cardinality reasons.
\end{proof}

For $n \in \N^+$, let $\zeta_n$  be the tight contact structures on $T^3$ defined as 
\[ \zeta_n = \ker (\sin(2 \pi n z)dx+ \cos (2 \pi n z)dy). \] 
\begin{theorem}[Giroux, \cite{giroux:2}] For any $n \in \N^+$, the contact  
structure $\zeta_n$ is universally tight and weakly symplectically fillable. 
Moreover $(T^3,\, \zeta_n)$ is contactomorphic to $(T^3,\, \zeta_m)$ if and only if $n=m$.
\end{theorem}
\begin{theorem}[\cite{kanda:1}, Theorem 0.1] Any tight contact structure $\xi$
on $T^3$ is contactomorphic to $\zeta_n$ for some $n$.
\end{theorem}
\begin{cor}\label{riempibili:T3}
All tight contact structures on $T^3$ are universally tight and weakly 
symplectically fillable.
\end{cor}

Take a primitive 
vector $(c_1, c_2, c_3) \in \Z^3$ with $c_3 \neq 0$ and complete it as the third row of 
a matrix 
$\Phi \in SL(3, \Z)$. The isotopy class of $\Phi^{-1}_* \zeta_n$ does not depend on the choice 
of the first and second rows of $\Phi$ because the stabiliser of $\zeta_n$ in 
$SL(3, \Z)$ acts transitively on them: see \cite{kanda:1} Theorem 0.2.

\begin{dfn}
Let $(c_1,c_2,c_3) \in \Z^3$ be a primitive vector and let $n$ be a
positive natural number. We set $\xi_{(n,c_1,c_2,c_3)} = \Phi^*
\zeta_n$.
\end{dfn}
  By \cite{kanda:1}, Theorem 7.6, 
$t(\xi_{(n,c_1,c_2,c_3)})=-|nc_3|$.
\begin{theorem}[\cite{honda:2}, Lemma 2.6] The tight contact structures
$\xi_{(n,c_1,c_2,c_3)}$ and $\xi_{(n',c_1',c_2',c_3')}$ are isotopic if and only if
$n=n'$ and $(c_1,c_2,c_3)= \pm (c_1',c_2',c_3')$. Moreover, any tight contact structure 
$\xi$ on $T^3$ with $t(\xi)<0$ is isotopic to $\xi_{(n,c_1,c_2,c_3)}$ for some 
$(n,c_1,c_2,c_3)$ with $c_3 \neq 0$.
\end{theorem}
\begin{theorem}[\cite{honda:2}, Section 2.5 Case 5]
If $e_0>0$ there is no tight contact structure $\xi$ on $T(e_0)$ with $t(\xi)<0$.
\end{theorem}
%%%%%%%%%%%%%%%%%%%%%%%%%%%%%%%%%%%%%%%%%%%%%%%%%%%%%%%%%%%%%%%%%%%%%%%%%%%%%%%
\subsection{Tight contact structures with $t=0$}
\begin{theorem}[\cite{honda:2}, Section 2.2 and Lemma 2.5]
The universally tight contact structures on $T(e_0)$ with maximal twisting number 
$t=0$ are in bijection with the set  $\N^+ \times \mathbb{P}(H_1(T^2; \Q))$.
\end{theorem}
The bijection in the theorem is given in the following way. $T(e_0)$ is also a 
$T^2$--bundle over $S^1$. Consider a convex $T^2$--fibre with infinite slope (i.e.
whose dividing curves are isotopic in $T(e_0)$ to $S^1$--fibres) and cut $T(e_0)$ 
along it
obtaining a $T^2 \times I$ with infinite boundary slopes. Make the boundary of
$T^2 \times I$ standard with horizontal ruling, and take a convex horizontal annulus
$A \subset T^2 \times I$. Gluing the boundary components of $A$ together, we obtain a torus
$T$ with a multicurve $\Gamma_T$. Let $n= \hbox{div}(T) \in \N^+$, and $s \in 
\mathbb{P}(H_1(T^2; \Q))$ the class of a connected component of $\Gamma_T$, then 
$(n,s)$ is the element in $\N^+ \times \mathbb{P}(H_1(T^2; \Q))$ associated to the 
tight contact structure on $T(e_0)$.

\begin{theorem}[\cite{ding-geiges:1}, Proposition 16] \label{riempibili:t=0} 
The universally tight contact structures on $T(e_0)$ with maximal twisting
number $t=0$ are all weakly symplectically fillable.
\end{theorem}

\begin{rem} When $e_0=0$, i.~e. when $T(e_0)=T^3$, the maximal twisting number $t$ 
reflects no geometric property of the tight contact structure, but depends only 
on the choice of a bundle structure $T^3 \to T^2$. 
\end{rem}

For $T^3$ we have
\[ \hbox{Tight}(T^3) \cong \N^+ \times \mathbb{P}(H_2(T^3; \Q)). \]
A tight $(T^3, \, \xi)$ corresponds to $(n, [T])$ such that $\xi$ is contactomorphic 
to $\zeta_n$ and $[T]$ is the unique homology class represented by a pre-Lagrangian 
torus in $(T^3, \, \xi)$. The set $\N^+ \times \mathbb{P}(H_1(T^2; \Q))$ of the isotopy 
classes of the tight contact structures on $T^3$ with maximal twisting number
$t=0$ embeds into $\N^+ \times \mathbb{P}(H_2(T^3; \Q))$ as $\N^+ \times H$, where 
$H \subset \mathbb{P}(H_2(T^3; \Q))$  is the hyperplane of the homology classes
represented by the fibred tori. 

\begin{theorem}[\cite{honda:2}, Proposition 2.3] There exist virtually 
overtwisted contact structures with $t=0$ on $T(e_0)$  only when $e_0>1$.  
There is one if $e_0=2$ and two if $e_0>2$.
\end{theorem}

 The virtually overtwisted contact structures with maximal twisting number
$t=0$ become overtwisted when pulled back to any covering of $T(e_0)$ induced by
a covering of the base $T^2$ and, by \cite{lisca-stipsicz:1}, are not weakly 
symplectically fillable.
%%%%%%%%%%%%%%%%%%%%%%%%%%%%%%%%%%%%%%%%%%%%%%%%%%%%%%%%%%%%%%%%%%%%%%%%%%%%%%
%%%%%%%%%%%%%%%%%%%%%%%%%%%%%%%%%%%%%%%%%%%%%%%%%%%%%%%%%%%%%%%%%%%%%%%%%%%
\section{Construction of the tight contact structures on $M(e_0,r)$}
\subsection{Thickening the singular fibre}

\begin{lemma} \label{decomposizione}
If $\xi$ is a tight contact structure on $M(e_0,r)$ with  maximal 
twisting number $t(\xi)$, then there exists a neighbourhood $V$ of the singular 
fibre $F$ such that $- \partial(M(e_0,r) \setminus V)$ is convex with slope $\frac{1}{t(\xi)}$.
 Moreover:
\begin{enumerate}
\item If $e_0 <0$, then $t(\xi) \geq -1$.
\item If $e_0 =0$, then $t(\xi)> - \frac{1}{r}$.
\item If $e_0 >0$, then $t(\xi)=0$.
\end{enumerate}
\end{lemma} 
\begin{figure}[ht!] 
\centering
\psfrag{U}{\footnotesize $U$}
\psfrag{V}{\footnotesize $V$}
\includegraphics[width=5.5cm]{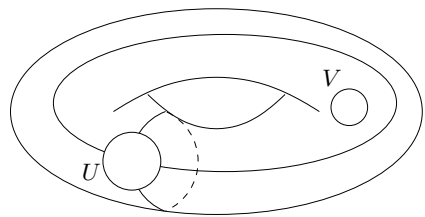}
\caption{How to cut $M \setminus (V \cup U)$}
\label{fig1.fig}
\end{figure}

\begin{figure}[ht!] 
\centering
\psfrag{V}{\footnotesize $V$}
\includegraphics[width=3.5cm]{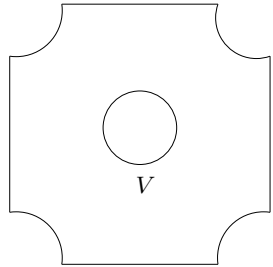}
\caption{The thickened torus with corners}
\label{fig2.fig}
\end{figure}

\begin{proof} 
In the following, we will call $M=M(e_0,r)$.
After an isotopy, we can find a Legendrian regular fibre $R$ with twisting number
 $t(\xi)$. The singular fibre $F$ can be made Legendrian with a very low twisting 
number $n$. We choose a standard neighbourhood $V$ of $F$ such that 
$- \partial (M \setminus V)$ has slope 
\[s_V = \frac{-n \beta - \beta'}{n \alpha + \alpha'} =- \frac{\beta}{\alpha} + \frac{1}{\alpha (n \alpha + \alpha')} < 
- \frac{\beta}{\alpha}\] 
where $\frac{\beta}{\alpha}=r$ and $\alpha'$, $\beta'$ are defined by $0 \leq \alpha' < \alpha$ and 
$\alpha' \beta - \alpha \beta'=1$.

If $t(\xi)=0$, choose a convex annulus $A$ so that one boundary component
 is a Legendrian ruling curve of $\partial (M \setminus V)$ and the other one is the 
Legendrian fibre $R$ with twisting number $t(\xi)$. By the imbalance principle 
\cite{honda:1} Proposition 3.17, we can perturb $A$ so that it contains a 
bypass attached to $\partial V$. By using this bypass we can thicken $V$ as far as 
there are singular 
points on $\partial A$, therefore we eventually get a solid torus $V$ 
with infinite boundary slope.

When $t(\xi)<0$, we choose a standard neighbourhood $U$ of $R$ such that $- \partial (M 
\setminus U)$ has boundary slope $s_U=-e_0+ \frac{1}{t(\xi)}$. 
In the convex annuli in figure \ref{fig1.fig}, whose boundary components are  
Legendrian ruling curves of $\partial (M \setminus U)$,  
all the dividing curves go from one boundary component to the other one, 
otherwise there would be 
 a bypass attached vertically to $U$ which would increase the twisting number 
of $R$ by the twisting number lemma.
When we cut $M \setminus (U \cup V)$ open along these two annuli, we obtain a
thickened torus with corners as shown in figure \ref{fig2.fig}. 

From slope $e_0 - \frac{1}{t(\xi)}$ on
$\partial (M \setminus U)$ by \cite{honda:1}, Lemma 3.11 we get, after 
rounding the edges, slope $e_0 +  
\frac{1}{t(\xi)}$, so the thickened torus we have obtained has boundary slopes 
$s_0= s_V <- r$ and $s_1= e_0 + \frac{1}{t(\xi)}$. If $e_0+ \frac{1}{t(\xi)} >  - r$, we have $s_1>s_0$ and there is an 
intermediate torus with infinite slope by \cite{honda:1}, Proposition 4.16. 
This torus would contradict the assumption about the maximality of the 
twisting number $t(\xi)$ of $R$, therefore $e_0+ \frac{1}{t(\xi)} \leq - r$. 
This implies that if $e_0>0$ than $t(\xi)=0$.
If $e_0+ \frac{1}{t(\xi)}=-r$, there is an overtwisted disc in a tubular 
neighbourhood of  the singular 
fibre with boundary on a Legendrian divide with slope $-r$.
We now divide into cases according to the sign of $e_0$. 
\begin{enumerate}

\item If $e_0<0$, then $e_0+ \frac{1}{t(\xi)} < -1 < - r$, therefore there is 
always an intermediate torus with slope $-1$ which forces the maximal twisting 
number $t(\xi)$  to be greater than or equal to $-1$. 

\item If $e_0=0$, then $\frac{1}{t(\xi)} < - r$. 
\end{enumerate}

In cases $1$ and $2$ we can find an intermediate convex torus with slope 
$\frac{1}{t(\xi)}$ because $\frac{1}{t(\xi)} \in [e_0+ \frac{1}{t(\xi)}, - r)$, and 
this torus bounds a neighbourhood $V$ of the singular fibre $F$ such  that 
$- \partial (M \setminus V)$ has slope $\frac{1}{t(\xi)}$.
\end{proof}

\begin{dfn}
Let $(M(e_0,r), \, \xi)$ be a tight contact manifold with maximal twisting number
$t(\xi)$, and let $V$ be a tubular neighbourhood of the singular fibre $F$
as in Lemma \ref{decomposizione} 
such that $- \partial (M \setminus V)$ has slope $\frac{1}{t(\xi)}$. Then the 
contact manifold $(M(e_0,r) \setminus V, \, \xi |_{M(e_0,r) \setminus V})$ will be called a 
{\em background} of $(M(e_0,r), \, \xi)$.
\end{dfn}   

\begin{dfn} If $\xi_0$ is a contact structure on $M \setminus V$ and $\eta$ is a contact 
structure on $V$ which match along the boundary, we will denote the glued
contact structure on $M$ by $\xi_0(\eta)$.
\end{dfn}

Generally, on a manifold with nonempty boundary we consider tight contact 
structures up to isotopies fixed on the boundary. On the contrary, in the 
classification of the backgrounds we will allow isotopies to move the boundary 
 because of the following lemma.
\begin{lemma}
Suppose that $\xi_1$ and $\xi_2$ are tight contact structures on $M$ and $V \subset M$ is a 
solid torus with convex boundary with respect to both $\xi_1$ and $\xi_2$. If 
$\xi_1 |_{M \setminus V}$ is isotopic to $\xi_2 |_{M \setminus V}$ by an isotopy not necessarily fixed at
the boundary, and $\xi_1 |_V$
is isotopic to $\xi_2 |_V$, then the contact structures $\xi_1$ and $\xi_2$ are isotopic.
\end{lemma}
\begin{proof}
Let $\phi_s$ be the isotopy of $M \setminus V$ such that $\phi_0$ is the identity and 
$(\phi_1)_* (\xi_1 |_{M \setminus V})= \xi_2 |_{M \setminus V} $. We can extend $\phi_s$ to $\widetilde{\phi}_s$ on all 
of $M$ so that $\widetilde{\phi}_0$ is the identity on $M$ and consider 
$(\widetilde{\phi}_1)_*(\xi_1)$.
By construction, $(\widetilde{\phi}_1)_*(\xi_1 |_{M \setminus V})= \xi_2 |_{M \setminus V}$, and by the 
classification of tight contact structures in \cite{honda:1},
$(\widetilde{\phi}_1)_*(\xi_1 |_V)$ is 
isotopic relative to the boundary to $\xi_2 |_V$ because they have the same
boundary slope and the same relative Euler class. Let $\psi_s$ be an isotopy 
between 
them, and $\widetilde{\psi}_s$ its extension to $M$ by putting it constantly equal 
to the identity outside $V$, then $\widetilde{\phi}_s \circ \widetilde{\psi}_s$ is an isotopy
between $\xi_1$ and $\xi_2$.
\end{proof}

%%%%%%%%%%%%%%%%%%%%%%%%%%%%%%%%%%%%%%%%%%%%%%%%%%%%%%%%%%%%%%%%%%%%%%%%%%%%%%%
\subsection{Tight contact structures with $t<0$}
In this section we present all tight contact structures $\xi$ on $M(e_0,r)$ 
with $t(\xi)<0$ as negative contact 
surgery on  fillable contact structures on $T(e_0)$. This result is obtained by
showing that the background of $(M(e_0,r), \xi)$ is contactomorphic to the 
complement of a standard neighbourhood of a vertical Legendrian curve in 
$T(e_0)$. 
For conciseness of notation, in the following we will often write $M$ instead
of $M(e_0,r)$.
\begin{prop}\label{background:e=0}
The background $(M \setminus V, \, \xi |_{M \setminus V})$ of $(M, \, \xi)$ with 
maximal twisting number $t(\xi)<0$ and integer Euler number $e_0=0$ is 
 contactomorphic to the complement  of a 
standard neighbourhood of a vertical Legendrian curve with twisting number
$t(\xi)$ in $(T^3, \, \xi_{(n,c_1,c_2,c_3)})$ for some $(n,c_1,c_2,c_3) \in \N^+
\times {\mathbb P}(H_2(T^3, 
\Q))$. Moreover, $(n,c_1,c_2,c_3)$ is uniquely determined by the dividing sets of
two non-isotopic, incompressible standard tori  intersecting along
 a common vertical Legendrian ruling curve with twisting number $t(\xi)$.
\end{prop}
\begin{proof}
We choose a vertical Legendrian curve $R$ with twisting number $t(\xi)$ in 
$M \setminus V$, and two standard tori $T_1$ and $T_2$ intersecting along $R$ as in 
the statement.
Let $n_i$ be the division numbers and let $\frac{p_i}{q_i}$ be the slope of $T_i$. 
These numbers satisfy the relations $-n_i q_i =t(\xi)$ for $i=1,2$ because 
$tb(R)=- \frac{1}{2} |R \cap \Gamma_{T_i} |$. 

Take a small standard neighbourhood $U$ of $R$ such that $T_i \cap \partial U$ is 
Legendrian. After cutting $(M \setminus V \cup U)$ along the two annuli $T_i \setminus U$ and 
rounding the edges as shown in Figures \ref{fig1.fig} and \ref{fig2.fig},
by \cite{honda:1}, Lemma 3.11  
we obtain a thickened torus $T^2 \times I$ with minimal boundary and 
 boundary slopes $\frac{1}{t(\xi)}$. This thickened torus is nonrotative,
otherwise an intermediate standard torus with slope $-r$ would produce an
overtwisted disc. By \cite{honda:1}, Lemma 5.7,
up to an isotopy which fixes one boundary
component, there is a unique 
nonrotative tight contact structure on $T^2 \times I$ with minimal boundary
and boundary slopes $\frac{1}{t(\xi)}$ , 
therefore there is at most one tight contact structure on $M \setminus V$ which induces
on $T_i$ a dividing set with division number $n_i$ and slope $\frac{p_i}{q_i}$ for 
$i=1,2$.

 Let $n=(n_1, n_2)$ be the greatest common divisor  
and set $c_1= - \frac{n_1 p_1}{n}$, $c_2= - \frac{n_2 p_2}{n}$ and $c_3= -\frac{t(\xi)}{n}$.
As their greatest common divisor is $(c_1, c_2, c_3)=1$, we can complete $\left (
\begin{matrix}
c_1 & c_2 & c_3 
\end{matrix} \right )$ to a matrix
\[\Phi = \left (
\begin{matrix}
a_1 & a_2 & a_3 \\
b_1 & b_2 & b_3 \\
c_1 & c_2 & c_3 \\
\end{matrix}
\right ) \in SL(3, \Z)\]
Fix coordinates $(x,y,z)$ on $T^3$, and consider the contact structure 
$\xi_{(n,c_1,c_2,c_3)}= \Phi^{-1}_* \zeta_n$.
 We claim that  $(M \setminus V, \, \xi |_{M \setminus V})$ is contactomorphic to 
 the complement of a standard neighbourhood of a vertical Legendrian curve 
in $(T^3, \, \xi_{(n,c_1,c_2,c_3)})$.
 In order to prove the claim, it is enough to show that the linear torus 
$T_1 \subset (T^3, \, \xi_{(n,c_1,c_2,c_3)})$ generated by 
$\left ( \begin{array}{c} 1 \\ 0 \\ 0 \end{array} \right )$ and 
$\left ( \begin{array}{c}  0 \\ 0 \\ 1 \end{array} \right )$ has division 
number $n_1$ and slope $\frac{p_1}{q_1}$, and the linear 
torus $T_2 \subset (T^3, \, \xi_{(n,c_1,c_2,c_3)})$ generated by 
$\left ( \begin{array}{c} 0 \\ 1
 \\ 0 \end{array} \right )$ and $\left ( \begin{array}{c}  0 \\ 0 \\ 1 
\end{array} \right)$ has division number $n_2$ and slope $\frac{p_2}{q_2}$. 
Equivalently, we can work with the tori $A(T_i) \subset (T^3, \zeta_n)$ generated by 
 $\left ( \begin{array}{c} a_i 
\\ b_i \\ c_i \end{array} \right )$ and $\left ( \begin{array}{c} a_3 \\ b_3 
\\ c_3 \end{array} \right )$ for $i=1,2$. Since $c_i \neq 0$ for $i=1,2,3$, there
is a  linear combination $X$ of 
$\frac{\partial}{\partial x}$  and $\frac{\partial}{\partial y}$ transverse to both $A(T_1)$ and $A(T_2)$. 
$X$ is a contact vector field of $(T^3, \zeta_n)$ for each $n$, and 
 the set $\Sigma = \{ p \in T^3 \; | \; X(p) \in \zeta_n(p) \}$ consists of  $2n$ parallel copies 
of a horizontal torus of the form $\{z \in \Z \}$. 

 The embeddings $\iota_i\co  T^2 \to T^3$ induced by the embeddings
$\widetilde{\iota}_i \co \R^2 \to \R^3$ given by
\[\widetilde{\iota}_i(u,v)= \left ( \begin{array}{c} 
ua_3+va_i \\ ub_3+vb_i \\ uc_3+vc_i
\end{array} \right )\]
are parametrisations of $A(T_i)$, for $i=1,2$.  
 The dividing set $\Gamma_{A(T_i)}= \Sigma \cap A(T_i)$ is the image of $2n$ parallel copies 
of the set $$\{ vc_i+uc_3 \in \Z \} = \bigcup \limits_{j=0}^{n_i \slash n} 
\{ -vp_i+uq_i \in \frac{jn}{n_i} \Z \},$$ 
which in turn consists of 
 $\frac{n_i}{n}$ parallel  
copies of  a  curves with 
slope $\frac{p_i}{q_i}$, therefore the dividing set $\Gamma_{A(T_i)}$ is the same 
dividing set induced by $\xi|_{M \setminus V}$ on $T_i$. 
\end{proof}
\begin{theorem}
Any tight contact structure $\xi$  on 
$M(0,r)$  with $t(\xi) \in (- \frac{1}{r}, 0)$ is a negative 
contact surgery on a vertical Legendrian curve with twisting number $t(\xi)$ 
in $(T^3, \, \xi_{(n,c_1,c_2,c_3)})$ for some $(n,c_1,c_2,c_3)$.
Conversely, any contact structure $\xi_{(n,c_1,c_2,c_3)}(\eta)$ on $M(0,r)$ obtained 
by negative contact surgery on $(T^3, \, \xi_{(n,c_1,c_2,c_3)})$ is tight.
\end{theorem}
\begin{proof}
The first half of the theorem comes from the previous proposition and from
   $- \frac{1}{r} < t(\xi)$.
 All  contact structures $\xi_{(n,c_1,c_2,c_3)}(\eta)$ obtained by negative
contact surgery on $(T^3, \xi_{(n,c_1,c_2,c_3)})$ are tight
 because all tight 
contact structures on $T^3$ are weakly symplectically fillable by Corollary
\ref{riempibili:T3}. 
 \end{proof}

\begin{theorem}
Let $\xi_{(n,c_1,c_2,c_3)}(\eta)$ be a tight contact structure on $M(0,r)$ obtained 
by negative contact surgery on a vertical Legendrian curve in the tight contact
manifold $(T^3, \, \xi_{(n,c_1,c_2,c_3)})$.
Let $\pi^*\xi_{(n,c_1,c_2,c_3)}(\eta)$ be the contact structure on 
$M(0,r, \ldots ,r)$ obtained as pull-back of $\xi_{(n,c_1,c_2,c_3)}(\eta)$  with respect to
the finite covering 
\[\pi \co M(0,r, \ldots ,r) \to M(0,r)\]
 induced by a finite covering of $T^2$. Then $\pi^*\xi_{(n,c_1,c_2,c_3)}(\eta)$ is tight.
\end{theorem}
\begin{proof}
Let $\widetilde{\xi}_{(n,c_1,c_2,c_3)}$ be the pull-back of $\xi_{(n,c_1,c_2,c_3)}$ with 
respect to the the finite covering of $T^3$ induced by the finite covering of 
$T^2$. 
By construction, $\pi^*\xi_{(n,c_1,c_2,c_3)}(\eta)$ is the contact structure 
$\widetilde{\xi}_{(n,c_1,c_2,c_3)}(\eta, \ldots ,\eta)$, obtained by negative contact surgery along 
a finite number of fibres of $T^3$.  

The contact structure $\widetilde{\xi}_{(n,c_1,c_2,c_3)}$ is tight
because all tight contact structures on $T^3$ are universally tight,
so it is also weakly symplectically fillable by Corollary
\ref{riempibili:T3}. The contact manifold $(M(0,r, \ldots ,r), \,
\widetilde{\xi}_{(n,c_1,c_2,c_3)}(\eta, \ldots ,\eta))$ is obtained by
negative contact surgery on a weakly symplectically fillable contact
manifold, therefore it is tight.
\end{proof} 

\begin{prop}
Let $\xi$ be a tight contact structure with maximal twisting number $t(\xi)=-1$ 
on the Seifert manifold $M=M(e_0,r)$ with integer Euler number $e_0 <0$. 
Then any background $(M \setminus V, \, \xi|_{M \setminus V})$ is contactomorphic to the complement of
 a standard neighbourhood of a vertical Legendrian curve
with twisting number $-1$ in $(T(e_0), \, \xi_0)$, where $\xi_0$ is a tight contact 
structure with $t(\xi_0)=-1$.
\end{prop} 
\begin{proof}
Let $T \subset M \setminus V $ be a standard vertical torus so that the manifold 
$M \setminus (T \cup V)$ is diffeomorphic to $\Sigma_0 \times S^1$, where $\Sigma_0$ is a pair of pants.
 We can assume that $T$ has vertical Legendrian ruling and its dividing set  
intersects the Legendrian ruling curves in two points. If this were not the case,
an annulus $A$ between a Legendrian ruling curve of $T$ and a Legendrian ruling 
curve  of $\partial (M \setminus V)$ would give a bypass along $T$ by the Imbalance Principle 
\cite{honda:1} Proposition 3.17, therefore we could decrease the number of 
intersection points between the dividing set and the Legendrian ruling curves of
$T$.

  Let $T_+$ and $T_-$ be the boundary tori of $\partial (M \setminus (V \cup T))$ corresponding to
$T$. $\xi|_{M \setminus (V \cup T)}$ is a tight contact structure with boundary slopes $1$ on 
$\partial (M \setminus V)$, $n$ on $T_+$,
and $-n+e_0$ on $T_-$. Since the sum of the slopes is $1+e_0 \leq 0$ and 
there are no vertical Legendrian curves with twisting number $0$, by 
\cite{honda:2}, Lemma 5.1
case 4(b), there are $1-e_0$ tight contact structures on $\Sigma_0 \times
S^1$ with those boundary slopes. Such contact structures are constructed
by removing a standard neighbourhood of a vertical Legendrian curve with 
twisting number $-1$ from a minimally twisting  $T^2 \times I$ with 
boundary slopes $n-e_0$ and $n$. Note here the effect of the orientation 
reversing  identification $T_- \cong T^2 \times \{ 0 \}$ on the slope. We can 
also assume that the standard 
neighbourhood of the vertical Legendrian curve is removed from an invariant 
collar of the boundary.

To have $M$ back from $M \setminus T$, we glue $T_+$ to $-T_-$ by the map  
$A(e_0)= \left ( \begin{matrix}
1 & 0 \\
- e_0 & 1 
\end{matrix} \right )$, therefore, by comparing with the construction in 
\cite{honda:2}, section 2.5, case 9,
 $(M \setminus V,\, \xi_{M \setminus V})$ is the complement of a 
vertical Legendrian curve with twisting number $-1$ in a circle 
bundle over the torus with Euler class $e_0$ with a tight contact structure
with maximal twisting number $t=-1$.

Given any slope $\overline{s} \in [s(-T_-), s(T_+)]$, we can find a 
convex torus $T' \subset M \setminus (T \cup V)$  with slope $\overline{s}$ such that $T_-$ and $T'$
 bound a thickened torus $T^2 \times [0, \frac{1}{2}] \subset M \setminus (T \cup V)$. 
Choose $\overline{s}= n-1$, then
remove $T^2 \times [0, \frac{1}{2}]$ from $M \setminus T$ and glue it back with $A(e_0)$ to 
the front, so that 
$M \setminus T'$  has 
boundary slopes  $n-e_0-1$ and $n-1$. Here one component of $\partial (M \setminus T')$ is
oriented with the outward normal and the other one with the inward normal.
In a similar way we can replace $n$ with 
$n+1$, so we have proved that the tight contact structure on $M$ does not 
depend on $n$.
\end{proof}

Conversely, given any tight contact structure $\xi_n$  on
$T(e_0)$ with $t(\xi_n)=-1$, for $n \in \Z / (1-e_0) \Z$ any negative contact surgery 
$(M(e_0,r), \, \xi_n(\eta))$ is tight. 

\begin{theorem} Let $e_0 < 0$. Any tight contact structure with $t=-1$ on 
$M(e_0,r)$  is 
negative contact surgery on a tight contact structure with $t=-1$ on $T(e_0)$.
Conversely, given any tight contact structure $\xi_n$  on
$T(e_0)$ with $t(\xi_n)=-1$, for $n \in \Z / (1-e_0) \Z$ any negative contact surgery 
$(M(e_0,r), \, \xi_n(\eta))$ is tight. 
\end{theorem}
\begin{proof}
Any tight contact structure with $t=-1$ on $M(e_0,r)$  is 
negative contact surgery on a tight contact structure $\xi_n$ with $t=-1$ on 
$T(e_0)$ because $- \frac{1}{r} < -1$. Conversely, any negative contact surgery
on $(T(e_0), \, \xi_n)$ is tight by \cite{ding-geiges:2}, Proposition 3 because 
$(T(e_0), \, \xi_n)$ is Stein fillable.
\end{proof}

\begin{theorem}
Let $\pi^*\xi_n(\eta)$ be the contact structure on $M(ke_0, r, \ldots ,r)$ obtained 
as pull-back of $\xi_n(\eta)$ with respect to a degree $k$ finite covering
\[\pi\co M(ke_0, r, \ldots ,r) \to M(e_0,r)\]
 induced by a covering of $T^2$. Then $\pi^*\xi_n(\eta)$ is tight.
\end{theorem}
\begin{proof}
By construction, $\pi^*{\xi_n(\eta)}= \widetilde{\xi}_n(\eta, \ldots ,\eta)$, where $\widetilde{\xi}_n$ is the
pull-back of $\xi_n$ to $T(ke_0)$. By \cite{honda:2}, Section 2.5, 
Case 9, $\widetilde{\xi}_n$ is a tight contact structure with maximal twisting number
$t(\widetilde{\xi}_n)=-1$, hence it is Stein fillable by Theorem \ref{riempibili}.
The contact manifold $(M(ke_0,r, \ldots ,r), \, \widetilde{\xi}_n(\eta, \ldots ,\eta))$ is tight 
because it is obtained by negative contact surgery on the Stein fillable 
contact manifold $(T(ke_0), \, \widetilde{\xi}_n)$.
\end{proof}
%%%%%%%%%%%%%%%%%%%%%%%%%%%%%%%%%%%%%%%%%%%%%%%%%%%%%%%%%%%%%%%%%%%%%%%%%%%%%%%
\subsection{Tight contact structures with $t=0$}
In this subsection we construct all tight contact structures $\xi$ on $M(e_0,r)$ 
with maximal twisting number $t(\xi)=0$. By Lemma \ref{decomposizione}, there is
 a tubular neighbourhood $V$ of the singular fibre such that 
$- \partial (M(e_0,r) \setminus V)$ is a convex torus with infinite slope. 
$M(e_0,r) \setminus V$ is diffeomorphic to $\Sigma \times S^1$, where $\Sigma$ is a punctured torus.
We will abusively identify $\Sigma$ with the image of a section $\Sigma \to \Sigma \times 
S^1$ and assume it is convex with Legendrian boundary and $\# \Gamma$--minimising in 
its isotopy class. 
 
The dividing set $\Gamma_{\Sigma}$ of $\Sigma$ consists of one arc with 
endpoints on $\partial \Sigma$ and 
some simple homotopically nontrivial closed curves.
\begin{dfn}
We define an {\em abstract dividing set} on an oriented surface $\Sigma$
as a multicurve $\Gamma_{\Sigma}$ together with a map $\pi_0(\Sigma \setminus
\Gamma_{\Sigma}) \to \{ +,- \}$ such that any connected component of $\Gamma_{\Sigma}$ belongs to the 
boundary of both a positive and a negative region. We say 
that an abstract dividing set is {\em tight} if its underlying multicurve does
 not have closed, homotopically trivial connected components. We say that it 
is {\em overtwisted} if it is not tight. 
\end{dfn}
In the following, we will almost always use the same symbol for both an 
abstract dividing set and for its underlying multicurve.  However, we will 
always specify what we are referring to, whenever it is relevant.

\begin{dfn}
Given an abstract dividing set $\Gamma_{\Sigma}$ on $\Sigma$, we denote by 
$\xi_{\Gamma_{\Sigma}}$ the $S^1$--invariant contact structure on 
$\Sigma \times S^1$ which induces the dividing set $\Gamma_{\Sigma}$ on a 
convex $\# \Gamma$--minimising section.
\end{dfn}

By Giroux's tightness criterion, \cite{honda:2} Lemma 4.2, $\xi_{\Gamma_{\Sigma}}$ is tight
(and in fact universally tight) if and only if $\Gamma_{\Sigma}$ is a tight abstract 
dividing set.
By \cite{honda:2}, Section 4.3, $(M(e_0,r) \setminus V, \, \xi |_{M(e_0,r) \setminus V})$ is 
contactomorphic to an $S^1$--invariant tight contact manifold 
$(\Sigma \times S^1, \, \xi_{\Gamma_{\Sigma}})$. We call $\eta = \xi|_V$ and   $\xi = \xi_{\Gamma_{\Sigma}}(\eta)$.

We recall that we have chosen the basis on $- \partial (M(e_0,r) \setminus V)$ so that $\partial \Sigma$ has slope
$e_0$ and the fibres have infinite slope and the basis on $- \partial (\Sigma \times S^1)$
so that $\partial \Sigma$ has slope $0$ and the fibres have infinite slope.

\begin{prop} \label{pluto} 
Let $\xi$ be  a tight 
contact structure on $M(e_0,r)$ with maximal twisting number $t(\xi)=0$ and fix 
a diffeomorphism $M \setminus V \cong \Sigma \times S^1$ so that $\Sigma$ is $\# \Gamma$--minimising. If $e_0 \leq 0$, 
then $\Gamma_{\Sigma}$ has no boundary parallel dividing curves. If $e_0 >0$ and $\Gamma_{\Sigma}$ has a
 boundary parallel dividing curve, then $\# \Gamma=1$.
\end{prop}

\begin{proof}
If $\Gamma$ contains  a boundary parallel dividing arc, then there is a singular 
bypass on $\Sigma $ by \cite{honda:1}, Proposition 3.18. By \cite{honda:1}, Lemma 
3.15, attaching this bypass to $- \partial (M \setminus V)$ we thicken $V$ 
to $V'$ so that $- \partial (M \setminus V')$ has slope $e_0$. 
 If $\# \Gamma \geq 2$, and 
$p \in \Sigma$
belongs to some other dividing curve, then $\{ p \} \times S^1$ is a Legendrian fibre 
with twisting number $0$ because $\xi|_{M \setminus V}$ is $S^1$--invariant by 
\cite{honda:2}, section 4.3.
 Applying the Imbalance principle, \cite{honda:1}, Proposition 3.17, we use 
this curve to find 
a vertical bypass attached to $\partial (M \setminus V')$. The attachment of this bypass 
gives a further thickening of $V'$ to $V''$ so that $- \partial (M \setminus V'')$ has 
infinite boundary slope again. By \cite{honda:1}, Proposition 4.16, there is 
a standard torus with slope $-r$ in $V'' \setminus V$. This torus produces an
overtwisted disc. 

If $\# \Gamma =1$, we pick a simple closed curve $C \subset \Sigma \setminus V'$ which 
does not disconnect $\Sigma$ and is disjoint from the dividing curve. By the 
Legendrian Realization Principle, \cite{honda:1}, Theorem 3.7, we  can 
arrange the characteristic foliation on $\Sigma$ so that $C$ is a closed leaf. 
Because of the $S^1$--invariance of $\xi |_{M \setminus V}$, $C \times S^1$ is a 
pre-Lagrangian torus with slope $0$. By Lemma \ref{prelagrangiano} we can 
perturb this torus in order to obtain a convex torus $T$ in standard 
form with slope $0$ and two dividing curves. The torus $T$ can be assumed to 
be disjoint from $V'$ because $C$ is disjoint from the boundary parallel 
dividing arc producing the bypass.

If we cut $M \setminus V'$ open along $T$, we obtain $\Sigma_0 \times S^1$, where $\Sigma_0$ is a 
pair of pants, and all the three boundary tori have slope $0$ calculated with
respect to the product structure on $\Sigma_0 \times S^1$. Let $T_{\pm}$ be
the two boundary tori corresponding to $T$, and take a convex vertical 
annulus $A$ with Legendrian boundary between $T_+$ and $T_-$. If the dividing 
curves 
on $A$ do not go from $T_+$ to $T_-$, then there is a vertical bypass along $T$.
The attachment of this bypass produces a torus $T'$ with infinite slope.
Using a vertical Legendrian divide of $T'$ we can thicken $V'$ to $V''$ so that
$- \partial (M \setminus V'')$ has infinite slope again, thus obtaining a standard 
torus with slope $-r$ in $V'' \setminus V$. Again, this torus produces an overtwisted 
disc. If the dividing curves on $A$ go
from one boundary component to the other, then, after cutting along $A$ and 
rounding the edges, by \cite{honda:1}, Lemma 3.11 we obtain a torus with slope
$-1$ parallel to $- \partial(\Sigma \times S^1)$
which has slope $e_0-1$ calculated with respect to the basis of 
$- \partial (M \setminus V)$.
If $e_0 \leq 0$, by \cite{honda:1}, Proposition 4.16, there is a convex torus 
with slope $-r$ parallel to $- \partial (M \setminus V)$ which gives an 
overtwisted disc.
\end{proof}
\begin{prop} \label{background:t=0}
Let $\Sigma$ be a punctured torus and $\Gamma_{\Sigma}$ a tight abstract 
dividing set on $\Sigma$ without  boundary parallel dividing arcs. Then
$(\Sigma \times S^1, \, \xi_{\Gamma_{\Sigma}})$ can be contact embedded into a tight contact manifold 
$(T(e_0), \, \xi_{\Gamma_{\Sigma}})$ as the complement of a standard neighbourhood of a 
vertical Legendrian curve with twisting number $0$.
\end{prop}
\begin{proof}
Take a curve $C \subset \Sigma$ so that $C$ intersect each dividing arc in one single 
point. If we make $C$ Legendrian using the Legendrian realisation principle
\cite{honda:1} Corollary 3.8, the torus $T= C \times S^1$ is in standard form 
with infinite slope because $\xi_{\Gamma_{\Sigma}}$ is $S^1$--invariant. The contact 
structure $\xi_{\Gamma_{\Sigma}}$ restricted to $\Sigma \times S^1 \setminus T$ is still 
$S^1$--invariant and $\Gamma_{\Sigma \setminus C}= \Gamma_{\Sigma} \setminus C$ is a 
$\# \Gamma$--minimising section of $\Sigma \times S^1 \setminus T$.
Let $S$ be the surface diffeomorphic to an annulus obtained by gluing a disc 
$D$ to the boundary component of $\Sigma \setminus C$  corresponding to $\partial \Sigma$, and let $\Gamma_S$ 
be the natural extension of $\Gamma_{\Sigma \setminus C}$ to an abstract dividing set on $S$. The 
$S^1$--invariant tight contact manifold $(S \times S^1, \, \xi_{\Gamma_S})$ is 
contactomorphic to an $I$--invariant tight contact structure on $T^2 \times I$ by 
\cite{honda:1}, Theorem 2.3(4) because $\Gamma_S$ consists of parallel arcs 
joining the two different boundary components of $S$. 
The $S^1$--invariant contact manifold 
$(D \times S^1, \, \xi_{\Gamma_D})$ is a tight solid torus with infinite boundary slope and 
$\# \Gamma_{\partial D \times S^1}=2$. By \cite{honda:1} Theorem 2.3 there is a unique tight 
contact structure with such boundary conditions on the solid torus, therefore 
it is contactomorphic to a standard neighbourhood of a Legendrian curve with 
twisting number $0$.
 Gluing $T^2 \times \{ 0 \}$ to $T^2 \times \{ 1 \}$ with the matrix 
$\left ( \begin{matrix}
1 & 0 \\
-e_0 & 1
\end{matrix} \right )$ 
we get a tight contact structure on $T(e_0)$ which we call $\xi_{\Gamma_{\Sigma}}$ again, 
then $(\Sigma \times S^1, \, \xi_{\Gamma_{\Sigma}})$ contact embeds in $(T(e_0), \, \xi_{\Gamma_{\Sigma}})$ as the 
complement of a vertical Legendrian curve with twisting number $0$.   
\end{proof}
\begin{theorem}\label{costruzione:t=0}
let $\Gamma_{\Sigma}$ be a tight abstract dividing set on a punctured torus $\Sigma$ without 
 boundary parallel dividing arcs, and let $\nu L \subset (T(e_0), \, \xi_{\Gamma_{\Sigma}})$ be a 
standard neighbourhood of a vertical Legendrian curve $L$ with twisting 
number $0$. Then, for any tight contact structure $\eta$ on $D^2 \times S^1$ whose 
characteristic foliation on $\partial D^2 \times \Sigma$ is mapped by 
\[A(r)\co \partial (D^2 \times S^1) \to - \partial (T(e_0) \setminus \nu L) \] 
to the characteristic 
foliation of $- \partial (T(e_0) \setminus \nu L)$, the contact structure $\xi_{\Gamma_{\Sigma}}(\eta)$ on 
$M(e_0,r)$ is tight.
\end{theorem}
\begin{proof} By Proposition \ref{background:t=0}, the contact manifold
$(M, \, \xi_{\Gamma_{\Sigma}}(\eta))$ is 
obtained by negative contact surgery on $(T(e_0), \, \xi_{\Gamma_{\Sigma}})$, which is a 
weakly symplectically fillable contact manifold by Theorem 
\ref{riempibili:t=0} because it is universally tight by $S^1$--invariance.
\end{proof}
\begin{prop}\label{ultima}
Let $\Gamma_{\Sigma}^+$ and $\Gamma_{\Sigma}^-$ be the two tight abstract dividing
sets on the punctured torus $\Sigma$ with underlying multicurve 
$\Gamma_{\Sigma}$ with no boundary parallel dividing arcs. Then, for any tight contact 
structure $\eta$ on $D^2 \times S^1$ as in Theorem \ref{costruzione:t=0}, 
$(M(e_0,r), \, \xi_{\Gamma_{\Sigma}^+}(\eta))$ is isotopic to $(M(e_0,r), \, \xi_{\Gamma_{\Sigma}^-}(\eta))$.
\end{prop}

Most of the proof of Proposition \ref{ultima} relies on the following lemma.
\begin{lemma}\label{ultimissimo}
Let $\Gamma_{T^2}$ be a tight abstract dividing set on $T^2$, and $\gamma_1$, $\gamma_2 \subset T^2$ 
dividing curves bounding a negative (positive) region $C \subset T^2$. Given points 
$p_i \in \gamma_i$, for $i=1,2$,  the curves $\{p_i \} \times S^1$ in $(T^3, \, \xi_{\Gamma_{T^2}})$ are 
Legendrian and have twisting number $tb(\{p_i \} \times S^1)=0$. If $L_1$ and
$L_2$ are positive (negative) stabilisations of $\{p_1 \} \times S^1$ and $\{p_2 \} \times S^1$ 
respectively, then they are contact isotopic.
\end{lemma}
\begin{proof}
The curves $\{p_i \} \times S^1$ are Legendrian because $(T^3, \, \xi_{\Gamma_{T^2}})$ is 
$S^1$--invariant.
 
Let $(T^2 \times [0, \frac 12 ], \, \xi)$ be a positive (negative) basic slice with 
standard boundary and boundary slopes $s_0=0$ and $s_{\frac 12 }= \infty$, contact 
embedded in
$(T^3, \, \xi_{\Gamma_{T^2}})$ so that $\{p_1 \} \times S^1$ is a Legendrian divide of 
$T^2 \times \{ \frac 12 \}$ and $T^2 \times \{ 0 \} \subset C \times S^1$. Make the Legendrian ruling of 
$T^2 \times \{ 0 \}$ vertical, and consider a convex vertical annulus $A$ between
$\{p_1 \} \times S^1 \subset T^2 \times \{ \frac 12 \}$, and a vertical Legendrian ruling curve of $T^2 \times \{ 0 \}$. 
The dividing set of $A$ 
consists of a single dividing arc with both endpoints on $T^2 \times \{ 0 \}$, 
and the simply connected region of $A \setminus \Gamma_A$ is positive (negative). Then, by 
\cite{etnyre:0}, Lemma 2.20, a vertical Legendrian ruling curve of $T^2 \times \{ 0 \}$
is a positive (negative) stabilisation of $\{ p_1 \} \times S^1$. From the 
well-definedness up to isotopy of the stabilisation,
it follows that $L_1$ is contact isotopic to a vertical Legendrian ruling curve
of $T^2 \times \{ 0 \}$. We can repeat the same argument with a basic slice 
<$(T^2 \times [0, \frac 12], \, \xi)$ with the same sign and the same boundary slopes 
so that  $\{p_2 \} \times S^1$ is a Legendrian divide of $T^2 \times \{ \frac 12 \}$, and 
conclude 
that $L_2$ is contact isotopic to a vertical Legendrian ruling of $T^2 \times \{ 0 \}$.
Then $L_1$ and  $L_2$ are Legendrian isotopic.
\end{proof}
 
\begin{proof}[Proof of Proposition \ref{ultima}] Let $V \subset M(e_0,r)$ be a 
tubular neighbourhood of the singular fibre such that 
$(M(e_0, r) \setminus V, \, \xi_{\Gamma_{\Sigma}^{\pm}}(\eta)|_{M(e_0, r) \setminus V})=(\Sigma \times S^1, \, \xi_{\Gamma_{\Sigma}^{\pm}})$, and 
$(V, \, \xi_{\Gamma_{\Sigma}^{\pm}}(\eta)|_V)$ is contactomorphic to $(D^2 \times S^1, \eta)$. We can find a 
solid torus $V' \subset V$ such that $- \partial (M(e_0,r) \setminus V')$ has slope $-1$ because 
$- \frac 1r < -1$. By proposition \ref{background:t=0}, 
$(\Sigma \times S^1, \, \xi_{\Gamma_{\Sigma}^{\pm}})$ is the complement of a vertical Legendrian curve 
$\{ p_1 \} \times S^1$ with twisting number $0$ with $p_1 \in \gamma_1$, where $\gamma_1$ is  the 
completion in $T^2$ of the dividing arc in $\Sigma$. Analogously, 
$(M(e_0, r) \setminus V', \, \xi_{\Gamma_{\Sigma}^{\pm}}(\eta)|_{M(e_0, r) \setminus V'})=(M(e_0, r) \setminus V', \, \xi_{\pm})$ is 
the complement of a vertical Legendrian curve $L_1$ with twisting 
number $-1$ which is a stabilisation of $\{ p_1 \} \times S^1$ by \cite{etnyre:0}, 
Lemma 2.20. The sign of the stabilisation is determined by the sign of the 
basic slice $V \setminus V'$. In order to fix the notation, let us suppose it is 
positive. We claim that $(M(e_0, r) \setminus V', \, \xi_+)$ is contact isotopic to 
$(M(e_0, r) \setminus V', \, \xi_{-})$. Since our argument will be semi-local, we can 
assume without loss of generality that $T(e_0)$ is a trivial $S^1$--bundle.

Let $\gamma_2$ be another dividing curve on $T^2$ such that $\gamma_1$ and $\gamma_2$ bound a 
positive region. Choose a point $p_2 \in \gamma_2$ and consider the vertical Legendrian 
curve $\{ p_2 \} \times S^1$ and its positive stabilisation $L_2$. Then the complement
of a standard neighbourhood of $L_2$ in $(T(e_0), \, \xi_{\Gamma_{\Sigma}^+})$ is 
contactomorphic to $(M(e_0, r) \setminus V', \, \xi_-)$. By Lemma 
\ref{ultimissimo}, $L_1$ and $L_2$ are Legendrian isotopic, therefore, by 
\cite{etnyre:0}, Theorem 2.12, there is a contact isotopy $\varphi_t \co (T(e_0), \, 
\xi_{\Gamma_{\Sigma}^+}) \to (T(e_0), \, \xi_{\Gamma_{\Sigma}^+})$ such that $\varphi_0= id$ and $\varphi_1(L_1)=L_2$. This 
implies that $\varphi_1(M(e_0,r) \setminus V')$ is the complement of a standard neighbourhood 
of $L_2$, so $(\varphi_t|_{M(e_0, r) \setminus V'})^* \xi_{\Gamma_{\Sigma}^+}$ is a $1$--parameter family  of 
contact structures on $M(e_0, r) \setminus V'$ all with the same boundary condition 
joining $\xi_+$ to $\xi_-$. By Gray's Theorem, this implies that 
$(M(e_0,r) \setminus V', \, \xi_+)$ and $(M(e_0,r) \setminus V', \, \xi_-)$ are contact isotopic. 
\end{proof}

\begin{theorem}\label{exceptionality} 
Let $\pi^*\xi_{\Gamma_{\Sigma}}(\eta)$ be the contact structure on $M(ke_0,r, \ldots ,r)$ 
obtained by pull-back of $\xi_{\Gamma_{\Sigma}}(\eta)$ with respect to a degree $k$ finite 
covering  
\[\pi \co M(ke_0,r \ldots ,r) \to M(e_0,r)\]
 induced by a covering of $T^2$. Then 
$\pi^*\xi_{\Gamma_{\Sigma}}(\eta)$ is tight when $\# \Gamma_{\Sigma} >1$ and  is
overtwisted when $\# \Gamma_{\Sigma}=1$.
\end{theorem}

\begin{proof}
By construction, $\pi^*\xi_{\Gamma_{\Sigma}}(\eta)= \xi_{\Gamma_{\widetilde{\Sigma}}}(\eta, \ldots ,\eta)$, where 
$\Gamma_{\widetilde{\Sigma}}$ is the pull-back of $\Gamma_{\Sigma}$ with respect to the finite cover of
$T^2$ (we remind that the inclusion $\iota\co \Sigma \to T^2$ can be lifted to an inclusion
$\widetilde{\iota}\co \widetilde{\Sigma} \to T^2$). If $\# \Gamma_{\Sigma} >1$, then $\Gamma_{\Sigma}$ contains no 
boundary parallel arcs. In this case,  $\Gamma_{\widetilde{\Sigma}}$ does not contain boundary parallel dividing curves either, then the contact manifold 
$(M(ke_0,r, \ldots ,r), \, \xi_{\Gamma_{\widetilde{\Sigma}}}(\eta, \ldots ,\eta))$ is
obtained by negative contact surgery on the weakly symplectically fillable 
contact manifold $(T(ke_0), \, \xi_{\Gamma_{\widetilde{\Sigma}}})$, therefore 
$\xi_{\Gamma_{\widetilde{\Sigma}}}(\eta, \ldots ,\eta))$ is tight. If $\# \Gamma_{\Sigma}=1$, then $\Gamma_{\widetilde{\Sigma}}$ consists 
of $k$ boundary parallel arcs, one for each boundary component of 
$\widetilde{\Sigma}$.
When $\# \Gamma >1$ a boundary parallel dividing arc produces an overtwisted disc
as in the proof of Proposition \ref{pluto}, therefore
$\xi_{\Gamma_{\widetilde{\Sigma}}}(\eta, \ldots ,\eta))$ is overtwisted. 
\end{proof}
%%%%%%%%%%%%%%%%%%%%%%%%%%%%%%%%%%%%%%%%%%%%%%%%%%%%%%%%%%%%%%%%%%%%%%%%%%%%%%%
\section{Classification of the generic tight contact structures}
\subsection{Tight contact structures with $t<0$}
\begin{theorem}
Let $e_0 <0$. The tight contact manifolds $(M(e_0,r), \, \xi_n(\eta))$ and 
$(M(e_0,r), \, \xi_m(\eta'))$ obtained by negative contact surgery on $(T(e_0), \, \xi_n)$
and $(T(e_0), \, \xi_m)$
respectively are isotopic if and only if $m=n$ and $\eta$ is isotopic to $\eta'$
relative to the boundary.
\end{theorem}
\begin{proof}
By Theorem \ref{ciccio} there are $|e_0-1|$ choices for the background and by 
Lemma \ref{solidi}, there are $|(d_0+1) \ldots (d_k+1)|$ choices for 
$\eta$. On the other hand, Theorem 5.4 in \cite{gompf:1} shows how to produce Stein
fillings $(W, J)$ for $M(e_0,r)$. As a smooth manifold, $W$ is the same for
all the Stein fillings, and choosing all possible rotation numbers in the Legendrian 
realisation of the surgery link presenting
$W$, by \cite{gompf:1}, Proposition 2.3, we obtain $|(e_0-1)(d_0+1) \ldots (d_0+k)|$ 
Stein structures on $W$ with different $c_1(J)$. By \cite{lisca-matic:2}, 
Corollary 4.2, these Stein structures induce $|(e_0-1)(d_0+1) \ldots (d_0+k)|$ 
mutually non isotopic tight tight contact structures on $T(e_0)$. 
\end{proof}
Now we turn our attention to the tight contact structures on $M(0,r)$ with
$t<0$. 
The result proved here
is a generalisation of Theorem 4.7 in \cite{gluing} for the part concerning 
the distinction of the tight contact structures obtained by negative contact 
surgery on $T^3$.
We start with a preliminary digression about negative contact surgery on 
nonrotative tight contact structures on $T^2 \times I$. We fix a $S^1$--bundle 
structure $T^2 \times  I \to S^1 \times I$ and, consequently, a Seifert fibration $M' \to S^1 \times 
I$ with one singular fibre $F$ on any manifold $M'$ obtained by surgery along 
a fibre of $T^2 \times I$.

\begin{prop} \label{principale}
Let $\xi$ be a nonrotative tight contact structure on $T^2 \times I$ and $r$ a rational
number such that $- \frac{1}{r} < t(\xi)$. Let $L \subset (T^2 \times I, \, \xi)$ be a vertical 
Legendrian curve with twisting number $tb(L)=t(\xi)$ and $(M', \, \xi(\eta))$ a contact 
manifold obtained from $(T^2 \times I, \, \xi)$ by contact surgery along $L$ with 
surgery coefficient $- \frac{1}{r}$ with respect to the canonical framing. Then 
$\xi(\eta)$ is tight and any two properly embedded convex vertical annuli $A_0$, 
$A_1$ with common 
Legendrian boundary are contact isotopic, possibly after perturbing the 
characteristic foliation of $A_1$. In particular, $\Gamma_{A_0}$ 
is isotopic to $\Gamma_{A_1}$ and $\xi(\eta)|_{M' \setminus A_i}$ is isotopic to $\eta$
for $i=0,1$. 
\end{prop}
\begin{proof}
The contact structure $\xi(\eta)$ is tight by \cite{ding-geiges:2}, Proposition 3 
because 
 $(T^2 \times I, \, \xi)$ can be contact embedded into a weakly symplectically fillable
contact manifold.

First we prove that, if $A_0$ and $A_1$ are disjoint from the 
surgery support, then $\Gamma_{A_0}$ is isotopic to $\Gamma_{A_1}$.
 In this case, we can think of $A_0$ and $A_1$ also as convex 
annuli in $(T^2 \times I,\, \xi)$. Passing to a finite covering in the horizontal 
direction, we can assume that $\xi$ is nonrotative with integer boundary slopes,
therefore, by \cite{honda:1}, Lemma 5.7, $\Gamma_{A_0}$ is isotopic to 
$ \Gamma_{A_1}$. 

 Now we turn to the proof of the general case. We can assume without loss of 
generality 
that one of the two annuli, say $A_0$, is disjoint from the surgery support.
If this is not the case, we introduce a third annulus $A_2$ disjoint from the 
surgery support and isotope $A_0$ to $A_2$ first, and then isotope $A_2$ to 
$A_1$.
By Isotopy discretisation \cite{gluing}, Lemma 3.10, there 
is a sequence of convex vertical annuli $A_0, \ldots ,A_{\frac{i}{n}}, 
\ldots ,A_{\frac{n}{n}}=A_1$
all with the same Legendrian boundary such that  $A_{\frac{i+1}{n}}$ 
is obtained from  $A_{\frac{i}{n}}$ by attaching a single bypass. 
By the following Lemma \ref{unaltro}, if $A_{\frac{i}{n}}$ is disjoint from the 
surgery support, we 
can find another contact surgery presentation of $(M', \, \xi(\eta))$ with 
 surgery support disjoint from both $A_{\frac{i}{n}}$ and $A_{\frac{i+1}{n}}$. 
Once $A_{\frac{i}{n}}$ 
and $A_{\frac{i+1}{n}}$ are disjoint from the surgery support, we can conclude 
that $\Gamma_{A_{\frac{i}{n}}}$ is isotopic to $ \Gamma_{A_{\frac{i+1}{n}}}$,
so the bypass between  $A_{\frac{i}{n}}$ and $A_{\frac{i+1}{n}}$ is trivial.
By the triviality of trivial bypass attachments, \cite{gluing}, Lemma 2.10, 
$\xi(\eta)$ restricted to the layer between $A_{\frac{i}{n}}$ and $A_{\frac{i+1}{n}}$ is 
invariant, therefore $A_{\frac{i}{n}}$ and $A_{\frac{i+1}{n}}$ are contact isotopic, 
possibly after perturbing the characteristic foliation of $A_{\frac{i+1}{n}}$.
\end{proof}
\begin{lemma}\label{unaltro}
Let $A_0$ and $A_1$ be convex annuli in $M'$ as in the statement of Proposition 
\ref{principale} such that they intersect only at the boundary. Also assume 
that $A_0$ is disjoint from the surgery support $V$. Then we can find another 
contact surgery presentation of $(M', \, \xi(\eta))$ such that the surgery support is 
disjoint from both $A_0$ and $A_1$.
\end{lemma}
\begin{proof}
Let $N$ be the component of $M' \setminus (A_0 \cup A_1)$ homeomorphic to 
$D^2 \times S^1$. By the monotonicity of the slope, \cite{honda:1}, Proposition 4.16, 
$M' \setminus N$ has boundary slope  
$\frac{p}{q} \in [\frac{1}{t(\xi)}, -r)$ because $M' \setminus N \subset M' \setminus A_0$. 
Let $\gamma$ be a vertical Legendrian curve contained in  $\partial (M' \setminus N)$, then 
$q \leq \frac{1}{2} \# (\gamma \cap \Gamma_{\partial (M' \setminus N)})$. In particular, if we take $\gamma \subset A_0$, 
then $\frac{1}{2} \# (\gamma \cap \Gamma_{\partial (M' \setminus N)})= - t(\xi)$, so $q \leq -t(\xi)$.
This is possible only if 
$\frac{p}{q} = \frac{1}{t(\xi)}$, therefore 
there is a solid torus $V' \subset M' \setminus N$ with convex boundary with slope 
$- \frac{1}{t(\xi)}$ such that 
$(M' \setminus V', \, \xi(\eta)|_{M' \setminus V'})$ is contactomorphic to the 
complement of a standard neighbourhood of a vertical Legendrian curve $L'$ 
with twisting number $tb(L')= t(\xi)$ in $(T^2 \times I, \, \xi)$.
 In fact, we can identify $M' \setminus V'$ and $T^2 \times I \setminus \nu L$ 
so that $\partial V$ corresponds to $\partial (\nu L)$. 
Then, both $\xi(\eta)|_{M' \setminus V'}$ and $\xi|_{T^2 \times I \setminus \nu L}$ have the
same boundary  slopes, induce the same dividing set on $A_0$ and,
after cutting along $A_0$ and rounding the edges, yield a
nonrotative tight contact structure on $T^2 \times I$ with slope 
$\frac{1}{t(\xi)}$. By \cite{honda:1}, Lemma 5.7, there is only one such tight 
contact structure on $T^2 \times I$ up to contactomorphism, therefore $\xi(\eta)|_{M' \setminus V'}$ 
and $\xi|_{T^2 \times I \setminus \nu L}$ are contactomorphic.
\end{proof}
Let $T_1$, $T_2 \subset (M(e_0,r), \, \xi_{(n,c_1,c_2,c_3)}(\eta))$ be convex 
tori in the direction $(x,z)$ and $(y,z)$ disjoint from the surgery support 
$V$ such that their intersection is a common 
vertical Legendrian ruling curve $R$ with twisting number $t$: see Figure 
\ref{fig1.fig}. Let $s_i$ be the slope of $T_i$ and $n_i$ its division number. The 
fact that $R= T_1 \cap T_2$ is a common Legendrian ruling curve implies that the 
intersection between $R$ and $\Gamma_{T_i}$ is minimal for both $i=1,2$. Let 
$U \subset V$ be a solid torus such that $- \partial (M \setminus U)$ is convex and has slope 
$\frac{1}{d_0+1}$, where $d_0=[- \frac{1}{r}]$.

 \begin{prop} \label{fondamentale} Let $T' \subset M=M(e_0,r)$ be a standard
 torus isotopic to $T_1$ with vertical Legendrian ruling. Then:
\begin{enumerate}
\item $T'$ has slope $s_1$ and division number $n' \geq n_1$.
\item Any convex torus $T''$ intersecting $T'$ in a vertical Legendrian 
ruling curve of $T'$ is contact isotopic to $T_2$, possibly after perturbing 
its characteristic foliation. 
\end{enumerate} 
\end{prop}
\begin{proof}
To simplify the notation, in the proof we will fix 
$\xi = \xi_{(n,c_1,c_2,c_3)}(\eta)$.
By Isotopy Discretisation \cite{gluing} Lemma 3.10, 
there
is a sequence of convex tori $T_1=T_{(1)}, \ldots ,T_{(n)}=T'$ such that $T_{(i+1)}$ is 
obtained from $T_{(i)}$ by attaching a bypass. In particular, $T_{(i+1)}$ and $T_{(i)}$
bound $N_i$ diffeomorphic to $T^2 \times I$. We can assume inductively that $T_{(i)}$
 satisfies:
\begin{enumerate}
\item $s(T_{(i)})=s_1$ and $\hbox{div} (T_{(i)}) \geq n_1$.
\item There is a solid torus $U_i \subset M \setminus T_{(i)}$ isotopic in $M$ to $U$ such that 
$- \partial (M \setminus U_i)$ is convex with slope $\frac{1}{d_0+1}$.
\item $\xi |_{U_i}$ is isotopic to $\xi |_U$ and $\xi |_{M \setminus U_i}$ is isotopic to 
$\xi |_{M \setminus U}$.
\item There is a convex vertical annulus $A_i \subset M \setminus T_{(i)}$ with Legendrian 
boundary on $\partial (M \setminus T_{(i)})$ such that $A$ closes to a convex 
torus $\overline{A}_i \subset M$ contact isotopic to $T_2$. 
\end{enumerate}

We observe that the inductive hypotheses are satisfied for $T_1$ taking as $A_0$
the annulus obtained by cutting $T_2$ open along $R= T_1 \cap T_2$, and $U_0=U$. 
Assumptions 2 and 3 
imply that $(M \setminus T_i, \, \xi|_{M \setminus T_i})$ is negative contact surgery along a 
vertical Legendrian curve in $T^2 \times I$ with a nonrotative tight contact 
structure, therefore all vertical annuli as in assumption 4 have the same 
dividing set by Proposition \ref{principale}.

A priori there are three kinds of transitions from $T_{(i)}$ to $T_{(i+1)}$:
\begin{enumerate}
\item $\hbox{div}(T_{(i)})= \hbox{div}(T_{(i+1)})=1$ and $s(T_{(i+1)}) \neq s_1$
\item $s(T_{(i+1)}) = s_1$ and $\hbox{div}(T_{(i+1)})= \hbox{div}(T_{(i)})+1$
\item $s(T_{(i+1)}) = s_1$ and $\hbox{div}(T_{(i+1)})= \hbox{div}(T_{(i)})-1$
\end{enumerate}

\medskip
{\bf Case 1}\qua We will prove that there are no transitions which change the 
slope. Suppose by contradiction that $\hbox{div}(T_{(i)})= \hbox{div}(T_{(i+1)})=1$ 
and $s(T_{(i+1)})=s_1' \neq s_1$. We can assume either that the bypass is attached to 
$T_{(i)}$ from the front and $s_1'<s_1$, or that the bypass is attached to $T_{(i)}$ 
from the back and $s_1<s_1'$. We describe only the first possibility because 
the second one is symmetric. 

Attaching bypasses coming from a convex vertical annulus with Legendrian 
boundary $S \subset M \setminus N_i$ as long as the Imbalance Principle can be applied, we 
eventually obtain tori $T_{(i)}'$ and $T_{(i+1)}'$  bounding
$N_i' \supset N_i$. The tori $T_{(i)}'$ and $T_{(i+1)}'$ have either infinite 
slope or have slopes $s(T_{(i)}')= \frac{p}{q} > s(T_{(i+1)}')= \frac{p'}{q}$ 
and a convex vertical annulus with Legendrian boundary $S' \subset M \setminus 
N_i'$ between $T_{(i)}'$ and $T_{(i+1)}$ contains no more boundary parallel 
dividing arcs. In the first case we have a vertical Legendrian curve with 
twisting number $0$ in $M \setminus (T_{(i)} \cup A_i)$.  This is excluded by 
the classification of tight contact structures on solid tori because, by the
inductive hypothesis, there is no such curve either in $(U_i, \xi|_{U_i})$ or 
in $(M \setminus (T_{(i)} \cup A_i \cup U_i), \, \xi |_{M \setminus (T_{(i)} \cup A_i \cup U_i)})$. In the second case, after 
cutting along $S$ and rounding the edges, we get  slope $\frac{p-p'-1}{q} \geq 0$ 
on $\partial (M \setminus (N_i' \cup S)) \subset M \setminus (T_{(i)} \cup A_i)$. This is also impossible because 
$M \setminus (T_{(i)} \cup A_i)$ has meridional slope $-r<0$, and the existence of a torus 
with non negative slope contained in $M \setminus (T_{(i)} \cup A_i)$ would imply the 
existence of a vertical Legendrian curve with twisting number $0$.

\medskip
{\bf Case 2}\qua Now we consider transitions which increase the division number. 
The main point here is to show that the surgery support can be assumed to be 
disjoint from the transition. Take convex vertical annuli $A_i' \subset N_i$ and 
$A_i'' \subset M \setminus N_i$ with common Legendrian boundary and call $B_i=A_i' \cup A_i'' \subset M \setminus T_{(i)}$. 
By Proposition \ref{principale}, $A_i$ is contact isotopic to $B_i$
and $\xi|_{M \setminus (T_{(i)} \cup A_i)}$ is isotopic to $\xi|_{M \setminus (T_{(i)} \cup B_i)}$. If we set 
$A_{i+1}=A_i'' \cup A_i' \subset M \setminus T_{(i+1)}$, then $\overline{A}_i$ is contact isotopic to 
$\overline{A}_{i+1}$. The solid torus obtained by
rounding the edges of $M \setminus (T_{(i)} \cup B_i)$ has boundary slope 
$- \frac{1}{k_i}$ for some positive integer $k_i < \frac{1}{r}$ and the solid 
torus obtained by rounding the edges of $M \setminus (N_i \cup A_i'')$ has boundary slope
$- \frac{1}{k_i'} \in [- \frac{1}{k_i}, -r)$ because $M \setminus (T_{(i)} \cup B_i)$ has
meridional slope $-r$. In $M \setminus (N_i \cup A_i'')$ there is a solid torus $U_{i+1}$ such 
that $- \partial (M \setminus U_{i+1})$ is a convex torus with slope $\frac{1}{d_o+1}$ because
$\frac{1}{d_o+1} \in [- \frac{1}{k_i'}, -r)$. Applying Proposition \ref{comoda} 
 to $M \setminus (T_{(i)} \cup B_i)$, we conclude that $\xi|_{U_i}$ is isotopic to 
$\xi|_{U_{i+1}}$ because slope $\frac{1}{d_o+1}$ is a border between continued 
fraction blocks in $M \setminus (T_{(i)} \cup B_i)$, and no shuffling can occur between the 
signs of basic slices belonging to different continued fraction blocks. For 
the same reason $\xi|_{M \setminus (T_{(i)} \cup A_i \cup U_i)}$ is isotopic to 
$\xi|_{M \setminus (T_{(i)} \cup B_i \cup U_{i+1})}$, then we can conclude that $\xi|_{M \setminus U_i}$
is isotopic to $\xi|_{M \setminus U_{i+1}}$.  

\medskip
{\bf Case 3}\qua Transitions which decrease the division number can be handled 
in the same way as transitions which increase it, we need only to show that 
no transition can decrease $\hbox{div}(T_{(i)})$ below $n_1$. Suppose that, on the
contrary, $\hbox{div}(T_{(i)})=n_1$ and $\hbox{div}(T_{(i+1)})=n_1-1$ and
take convex vertical annuli $A_i' \subset N_i$ and $A_i'' \subset M \setminus N_i$ with common 
Legendrian boundary.  The dividing set of $B_i= A_i' \cup A_i'' \subset M \setminus T_{(i)}$ has at 
least one boundary parallel dividing arc, but the same total number of 
dividing arcs as $A_0$. This is a contradiction because, by Proposition 
\ref{principale} and the inductive hypothesis, the dividing set 
$\Gamma_{\overline{B}_i}$ on the torus $\overline{B}_i$ obtained by gluing the 
boundary components of $B_i$ is isotopic to $\Gamma_{\overline{A}_0}= \Gamma_{T_2}$ and 
$\Gamma_{A_0}$ contains no boundary parallel dividing arcs. 

It remains to prove that any convex torus $T''$ isotopic to $T_2$ and 
intersecting $T'$ in a Legendrian ruling curve of $T'$ is contact isotopic to 
$T_2$. The Legendrian curves $T' \cap \overline{A}_n$ and $T' \cap T''$ are Legendrian 
isotopic because they are both Legendrian ruling curves of $T'$. Let 
$\varphi_t \co M \to M$ be a contact isotopy extending 
the Legendrian isotopy between  $T' \cap T''$ and $T' \cap \overline{A}_n$, so that
$T' \cap \overline{A}_n = T' \cap \varphi_1 (T'')$. By Proposition \ref{principale}, $A_n=
\overline{A}_n \setminus T'$  is contact isotopic to $\varphi_1(T'') \setminus T'$, therefore the proof 
is finished because $\overline{A}_n$ is contact isotopic to $T_2$.
\end{proof}
\begin{theorem} Let $\xi_{(n,c_1,c_2,c_3)}$ and $\xi_{(n',c_1',c_2',c_3')}$ be
 tight contact structures over $T^3$ with $c_3 \neq 0$ and $c_3' \neq 0$. Then the tight
contact structures 
$\xi_{(n,c_1,c_2,c_3)}(\eta)$ and $\xi_{(n',c_1',c_2',c_3')}(\eta')$ over $M(0,r)$ constructed
by negative contact surgery are isotopic if and only if $n=n'$,
$(c_1,c_2,c_3)= \pm (c_1',c_2',c_3')$ and $\eta$ is isotopic to $\eta'$.
\end{theorem}
\begin{proof} Suppose $\xi_{(n,c_1,c_2,c_3)}(\eta)$ is isotopic to 
$\xi_{(n',c_1',c_2',c_3')}(\eta')$ and call it $\xi$. Because of the presentation of $\xi$ as
$\xi_{(n,c_1,c_2,c_3)}(\eta)$, in $M=M(0,r)$ we find  a solid torus $V \subset M$
so that $- \partial (M \setminus V)$ is convex with slope $\frac{1}{t}= - \frac{1}{|nc_3|}$, 
 $\xi|_V = \eta$ and $M \setminus V$ can be contact embedded in $(T^3, \, \xi_{(n,c_1,c_2,c_3)})$
as the complement of a vertical Legendrian curve. In $M \setminus V$ we choose
tori $T_1$ and $T_2$ with slope $s(T_i)= \frac{c_i}{c_3}$ and division number 
$\hbox{div}(T_i)= n(c_i,c_3)$ respectively intersecting along a common vertical 
Legendrian ruling curve. See the proof of Proposition \ref{background:e=0} for 
details. Similarly, because of 
the the presentation of $\xi$ as $\xi_{(n',c_1',c_2',c_3')}(\eta')$, in $M=M(0,r)$ we 
find tori $T_1'$ and $T_2'$ with slope $s(T_i')= \frac{c_i'}{c_3'}$ and division 
number $\hbox{div}(T_i')= n'(c_i',c_3')$ respectively and a solid torus 
$V' \subset M \setminus (T_1' \cup T_2')$ such that $- \partial (M \setminus V')$ is convex with slope 
$\frac{1}{t'}= - \frac{1}{|n'c_3'|}$ and $\xi|_V' = \eta'$.

Proposition \ref{fondamentale}, implies that $\Gamma_{T_i}$ is isotopic to 
$\Gamma_{T_i'}$ for $i=1,2$, therefore $n=n'$ and $(c_1, c_2, c_3)= \pm
(c_1', c_2', c_3')$. Applying Proposition \ref{fondamentale} to the first cut,
and Proposition \ref{principale} to the second one, we prove that $\eta \cong 
\xi|_{M \setminus (T_1 \cup T_2)}$ is isotopic to $\eta' \cong \xi|_{M 
\setminus (T_1' \cup T_2')}$, thus concluding the proof.
\end{proof}

%%%%%%%%%%%%%%%%%%%%%%%%%%%%%%%%%%%%%%%%%%%%%%%%%%%%%%%%%%%%%%%%%%%%%%%%%%%%%%%
\subsection{Tight contact structures with $t=0$} 
The aim of this section is the classification of the tight contact structures
 on
$M(e_0,r)$ constructed by negative contact surgery on a vertical Legendrian
curve in a tight, $S^1$--invariant contact structure on $T(e_0)$.
   
Given two multicurves $\Gamma$ and $\Gamma'$ on a surface $\Sigma$, we say that
they are {\em diffeomorphic} if there exists a diffeomorphism $\phi \co \Sigma 
\to \Sigma$ such that $\phi |_{\partial \Sigma}= id$, and $\phi (\Gamma)= 
\Gamma'$. If we consider $\Gamma$ and $\Gamma'$ as abstract dividing sets, we 
require
in addition that $\phi$ maps positive regions to positive regions and negative 
regions to negative regions. 

Let $\Sigma_0$ be a pair of pants with $\partial \Sigma_0 =C_0 \cup C_1 \cup C_2$, and let $\Gamma_{\Sigma_0}$ be 
an abstract dividing set on $\Sigma_0$ with $\# \Gamma_{\Sigma_0} \cap C_i \neq \emptyset$ for 
$i= 0,1,2$. If $\# \Gamma_{\Sigma_0} \cap C_2 =2$, there is 
a canonical way  up to isotopy to extend 
$\Gamma_{\Sigma_0}$ to an abstract dividing set in $A= S^1 \times I$, namely 
by gluing a disc $D$ to $C_2$ along
the boundary and joining the endpoints of $\Gamma_{\Sigma_0}$ on $\partial 
\Sigma_0$ with an arc 
contained in $D$. We will call the extension $\Gamma_A$, and will 
denote by $\xi_{\Gamma_A}$ the contact structure on $T^2 \times I$ which
is $S^1$--invariant over $\Gamma_A$. 
The contact manifold $(M', \, \xi_{\Gamma_{\Sigma_0}}(\eta))$ obtained by 
negative contact surgery on a vertical Legendrian curve with twisting number
$0$ in $(T^2 \times I, \, \xi_{\Gamma_A})$ is tight if and only if
$\Gamma_A$ is tight. In fact, if $\Gamma_A$ is a tight abstract dividing set,
we can take a tight abstract dividing set $\Gamma_{T^2}$ on $T^2$ together with an 
embedding $\iota \co A \to T^2$ so that $\iota (\Gamma_A) = \Gamma_{T^2} \cap \iota(A)$.
By Giroux Tightness Criterion \cite{honda:1}, Lemma 4.2, $(T^3, \, \xi_{\Gamma_{T^2}})$ is 
universally tight, therefore it is also symplectically fillable by Theorem 
\ref{riempibili:t=0}. Since $(T^2 \times I, \, \xi_{\Gamma_A})$ contact embed into 
$(T^3, \, \xi_{\Gamma_{T^2}})$, then $(M', \, \xi_{\Gamma_{\Sigma_0}}(\eta))$ is tight because it can be 
embedded into a contact manifold obtained by negative contact surgery on 
$(T^3, \, \xi_{\Gamma_{T^2}})$.

On the contrary, if $\Gamma_A$ is overtwisted, then $\Gamma_{\Sigma_0}$ 
contains either a homotopically trivial closed curve, or a boundary parallel 
dividing arcs with endpoints on $C_2$. Then $(M', \, \xi_{\Sigma_0}(\eta))$ is 
overtwisted by the same argument as in the proof of Proposition \ref{pluto} 
because $\# \Gamma_{\Sigma_0}>1$.

\begin{lemma}\label{confronto} 
Let $\Gamma_{\Sigma_0}$ and $\Gamma_{\Sigma_0}'$ be abstract dividing sets  on the pair of pants
$\Sigma_0$ such that their completions $\Gamma_A$ and $\Gamma_A'$ are 
tight, and let $\widehat{\Gamma}_A$, $\widehat{\Gamma}_A'$ be 
 obtained  from $\Gamma_A$, $\Gamma_A'$ by 
 throwing away every pair of closed curves bounding an annulus. 
If $\xi_{\Gamma_{\Sigma_0}}(\eta)$ is isotopic to $ \xi_{\Gamma_{\Sigma_0}'}(\eta')$,
then $\widehat{\Gamma}_A$ is diffeomorphic to  $\widehat{\Gamma}_A'$ 
(as abstract dividing sets) and $\eta$ is isotopic to $\eta'$.
\end{lemma}
\begin{proof}
Let $\xi_{\Gamma_{\Sigma_0}}(\eta)$ and $\xi_{\Gamma_{\Sigma_0}'}(\eta)$ be isotopic tight contact structures.
 By definition  there exist 
neighbourhoods $V$ and $V'$ of the singular fibre and sections $\sigma\co \Sigma_0 \to M' \setminus V$
and $\sigma'\co \Sigma_0 \to M' \setminus V'$ coinciding on a neighbourhood of $\partial \Sigma_0$ such that:
\begin{enumerate}
\item $- \partial (M \setminus V)$ and $- \partial (M \setminus V')$ have infinite boundary slope. 
\item $\xi_{\Gamma_{\Sigma_0}}(\eta)|_{V}$ is isotopic to $\eta$. 
\item $\xi_{\Gamma_{\Sigma_0}'}(\eta')|_{V'}$ is isotopic to  $\eta'$. 
\item $\sigma(\Sigma_0)$ and $\sigma'(\Sigma_0)$ are convex with Legendrian boundary,  
$\Gamma_{\sigma(\Sigma_0)}= \Gamma_{\Sigma_0}$ and $\Gamma_{\sigma'(\Sigma_0)}= \Gamma_{\Sigma_0}'$. 
\end{enumerate}
If we glue a $S^1$--invariant tight contact structure on 
$T^2 \times I$ to either component of $\partial M'$, the result is tight if 
and only if $\Gamma_A$ (respectively $\Gamma_A'$)
 glued to the dividing set on a convex horizontal annulus in $T^2 \times I$
produces no homotopically trivial curves.
We will exploit this fact to prove the lemma using a technique called
{\em Template Attaching}, first introduced by Honda in \cite{honda:1}, 
section 5.3.2. In the following we will call {\em elementary template} a 
thickened torus $T^2 \times I$ carrying a tight contact structure which is
 $S^1$--invariant over a horizontal annulus with only one boundary 
parallel dividing arc and $2(n-1)$ dividing arcs with endpoints on different 
boundary components.  

The set $\partial \Gamma_A= \partial \Gamma_A'$ consists of a finite collection of points with 
cardinality  
$2 (\# \widehat{\Gamma}_A)= 2 (\# \widehat{\Gamma}_A')$. Given two points 
$p,q \in \partial \Gamma_A$ (respectively $\partial \Gamma_A'$) joined by an arc in the dividing set, we 
denote by $(p,q)$ the arc in $\Gamma_A$ (respectively $\Gamma_A'$) joining them. We 
partition $\partial \Gamma_A$ in the two subsets $\partial \Gamma_A \cap C_0= \{ p_0, \ldots ,p_n \}$ and $\partial \Gamma_A \cap C_1 = 
\{ p_0', \ldots ,p_m' \}$ and put a cyclic order on them.

We work by induction on the number of dividing arcs in $\Gamma_A$ with both 
endpoints on the same boundary component. The base step is when there are no 
such curves, or when $\# \widehat{\Gamma}_A=2$. In the first case, $\Gamma_A$ coincides with 
$\widehat{\Gamma}_A$, and defines an order preserving bijection 
$\{ p_0, \ldots ,p_n \} \to \{ p_0', \ldots ,p_n' \}$, which 
determines it up to diffeomorphism. We claim that $\Gamma_A'$ has no boundary
parallel dividing arcs either and induces the same bijection as $\Gamma_A$.

If $ \widehat{\Gamma}_A$ contains no boundary parallel dividing arcs, then no single
elementary template attaching produces a homotopically trivial closed curve. 
Suppose $\Gamma_A'$ contains a boundary parallel arc $(p_i, p_{i+1})$. Then the 
attachment of an elementary template such that $p_i$ and $p_{i+1}$ are the 
endpoints of the boundary parallel dividing arc in its horizontal annulus
produces an overtwisted disc, giving a contradiction. See Figure 
\ref{template2.fig} case (b).

\begin{figure}[ht!] 
\centering
\psfrag{a}{\scriptsize (a)}
\psfrag{b}{\scriptsize (b)}
\psfrag{c}{\scriptsize (c)}
\includegraphics[width=12cm]{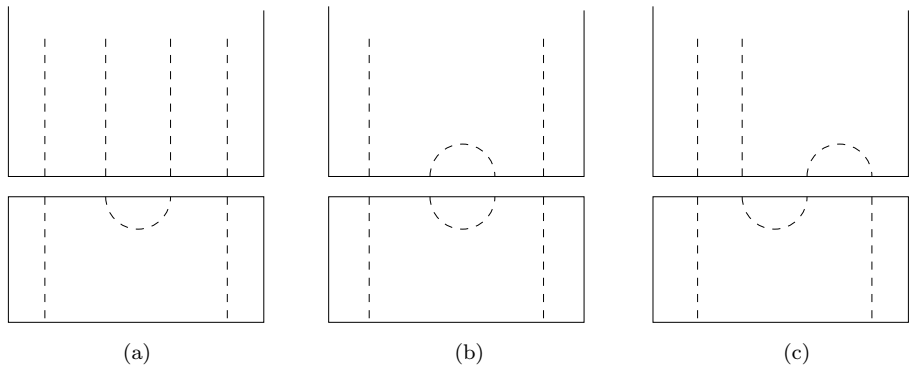}
\caption{The attaching of an elementary template:\qua Cases (a) and (c) preserve
tightness, cases (b) produces an overtwisted disc.}
\label{template2.fig}
\end{figure}

Let $T_0$ and $T_1$ be the two components of $\partial M'$. If we attach  elementary 
templates to $T_0$ and $T_1$ such that $\{ p_0, p_1 \}$ and $\{ p_i', p_{i+1}' \}$ are the 
endpoints of the boundary parallel dividing arcs, then the only case in which
 we get an overtwisted disc is when there are dividing arcs $(p_0, p_i')$ and 
$(p_1, p_{i+1}')$. This must be true for both $\Gamma_A$ and $\Gamma_A'$, therefore
the two dividing sets are isomorphic.

When $\# \widehat{\Gamma}_A=2$, we have to distinguish the cases when $\Gamma_A$ 
consists of two 
non boundary parallel dividing arcs, or when it consists of one boundary 
parallel dividing arc
on each side and a number of closed curves. In the first case, 
$(M', \, \xi_{\Gamma_{\Sigma_0}}(\eta))$ remains tight after gluing $S^1$--invariant tight contact 
structures with a boundary parallel dividing curve on its horizontal annulus 
in any possible  way. In the second case some gluing produce an overtwisted 
disc. This forces $\widehat{\Gamma}_A'$ to be diffeomorphic to 
$\widehat{\Gamma}_A$. We observe that 
template attaching cannot detect multiple closed curves in $\Gamma_A$ and 
$\Gamma_A'$. This is the reason why we work with $\widehat{\Gamma}$ instead of with $\Gamma$.

Now we suppose the lemma true when $\Gamma_A$ has $k-1$ arcs with endpoints on
the same boundary component. Let $\Gamma_A$ have $k$ of such arcs, and suppose 
that $(p_i, p_{i+1})$ is one of them. If we glue an elementary template to 
$(M', \, \xi_{\Gamma_{\Sigma_0}}(\eta))$ so that the boundary parallel dividing arc in its 
horizontal annulus matches with $(p_i, p_{i+1})$ to give a closed homotopically 
trivial curve, we produce an overtwisted disc. This fact implies that 
$(p_i, p_{i+1})$ is also a boundary parallel dividing arc in $\Gamma_A'$. 

After slightly perturbing $\xi_{\Gamma_{\Sigma_0}}(\eta)$ and $\xi_{\Gamma_{\Sigma_0}'}(\eta')$, for both 
contact structures $A$ contains a bypass along $\partial M'$ corresponding to the 
boundary parallel dividing arc $(p_i, p_{i+1})$. After attaching these bypasses to
 $\partial M'$, and removing the collars of $\partial M$ in which the bypass 
attachment takes place, we obtain manifolds $M_1'$, $M_2'$ with tight contact 
structures  $\xi_{\Gamma_{\Sigma_0}}(\eta))|_{M_1'} = \xi_{\Gamma_{\Sigma_0} \setminus (p_i,p_i+1)}(\eta)$ and 
$\xi_{\Gamma_{\Sigma_0}'}(\eta'))|_{M_2'}=\xi_{\Gamma_{\Sigma_0}'\setminus (p_i,p_i+1)}(\eta')$.
If we prove that $\xi_{\Gamma_{\Sigma_0} \setminus (p_i,p_i+1)}(\eta)$ is isotopic to 
$\xi_{\Gamma_{\Sigma_0}'\setminus (p_i,p_i+1)}(\eta')$ we can use the inductive hypothesis to conclude that
$\widehat{\Gamma}_A\setminus (p_i,p_i+1)$ is diffeomorphic to $\widehat{\Gamma}_A'\setminus (p_i,p_i+1)$. From this 
it follows that $\widehat{\Gamma}_A$ is diffeomorphic to $\widehat{\Gamma}_A'$.

To prove that $\xi_{\Gamma_{\Sigma_0} \setminus (p_i,p_i+1)}(\eta)$ is isotopic to 
$\xi_{\Gamma_{\Sigma_0}'\setminus (p_i,p_i+1)}(\eta')$, we glue an elementary 
template to $M'$ so that the boundary 
parallel dividing arc on its horizontal annulus joins $p_{i+1}$ to $p_{i+2}$. The 
resulting contact manifold $M''$ is tight and both $M'' \setminus M_1'$ and $M'' \setminus M_2'$ 
are contactomorphic to $I$--invariant thickened tori. See Figure 
\ref{template2.fig} case (c).

If $M$ decomposes as $(M \setminus V) \cup V$, consider the inclusions $\iota_V \co V \hookrightarrow M$ and 
$\iota_{M \setminus V}\co M \setminus V \hookrightarrow M$. From the obstruction theoretical definition of the Euler
class it is immediate that 
\[(\iota_V)_* PD (e(\xi|_V, s))+(\iota_{M \setminus V})_* PD(e(\xi|_{M \setminus V}, s)= PD e(\xi)\] 
for any section $s$ of $\xi$ on $\partial V$. 
If $\xi_{\Gamma_{\Sigma_0}}(\eta)$ and $\xi_{\Gamma_{\Sigma_0}}(\eta')$ are isotopic, then $e(\xi_{\Gamma_{\Sigma_0}}(\eta))=
e(\xi_{\Gamma_{\Sigma_0}}(\eta'))$. Moreover, $\xi_{\Gamma_{\Sigma_0}}(\eta)|_{M \setminus V}= \xi_{\Gamma_{\Sigma_0}}(\eta')|_{M \setminus V}= 
\xi_{\Gamma_{\Sigma_0}}$, then 
\[
(\iota_{M \setminus V})_* PD(e(\xi_{\Gamma_{\Sigma_0}}(\eta)|_{M \setminus V}, s))= (\iota_{M \setminus V})_* PD(e(\xi_{\Gamma_{\Sigma_0}}(\eta')|_{M \setminus V}, s)).
\] 
By difference, $(\iota_V)_* PD(\eta, s)= (\iota_V)_* PD(\eta', s)$,
therefore $e(\eta, s)=e(\eta', s)$ because $(\iota_V)_*$ is injective. By  \cite{honda:1}, 
Proposition 4.23, this proves that
$(D^2 \times S^1, \eta)$ and $(D^2 \times S^1, \eta)$ are isotopic
\end{proof}
 
Given $(M, \, \xi)$, and a neighbourhood $V$ of the singular fibre $F$ such that 
$- \partial (M \setminus V)$ is a standard torus with infinite slope, we can 
modify the Seifert fibration $\pi \co M \to T^2$ by an isotopy so that 
$M \setminus V$ fibres onto $T^2 \setminus D= \Sigma$, where $D$ is an embedded
disc. Let $\sigma \co \Sigma \to M \setminus V$ be a section such that
$\sigma ( \Sigma)$ is convex with Legendrian boundary and 
$\# \Gamma$--minimising in its isotopy class. With an abuse of notation, we
will denote $\sigma(\Sigma)$ simply by $\Sigma$, and its dividing set by $\Gamma_{\Sigma}$.
We will denote by $\Gamma$ the extension of $\Gamma_{\Sigma}$ to $T^2$ obtained by joining the 
endpoints of $\Gamma_{\Sigma}$ with an arc in $D$. 
\begin{prop}\label{conto}
Let $\gamma \subset T^2$ be a homotopically nontrivial simple closed curve  
disjoint from the image of the singular fibre, and let ${\mathcal T}_{\gamma}$ be the 
family of convex or pre-Lagrangian tori in $(M(e_0,r), \, \xi_{\Gamma_{\Sigma}}(\eta))$ isotopic to 
$\pi^{-1} (\gamma)$. If we define the division number of a pre-Lagrangian torus to be 
zero, then the equality
\[
\min \limits_{T \in {\mathcal T}_{\gamma}} \hbox{\rm div} (T)= \frac{1}{2} | \gamma \cap \Gamma |
\] 
holds.
\end{prop}

\begin{rem} We say that a convex vertical torus $T$ in $M(e_0,r)$ 
has infinite slope if its dividing set is isotopic in $M(e_0,r)$ to regular 
fibres, otherwise we say that $T$ has finite slope. In general we cannot give 
a well-defined value to the slope of $T$ when it is finite.
\end{rem}

\begin{lemma} 
${\mathcal T}_{\gamma}$ 
contains a pre-Lagrangian torus if and only if it contains a convex torus which
does not have infinite slope.
\end{lemma}
\begin{proof}
Suppose there is a pre-Lagrangian torus $T \in {\mathcal T}_{\gamma}$: then after a 
suitable choice of coordinates $(x,y,z)$ in a neighbourhood of $T$ so that
$T= \{ y=0 \}$, the contact structure has equation
$dz -y \, dx$ in a neighbourhood of $T$. This local model shows that for some
small $\epsilon \neq 0$, the torus $T_{\epsilon}$ is pre-Lagrangian and has 
rational slope different from the slope of $T$.
Then we perturb $T_{\epsilon}$ and obtain a convex torus $T'$ with finite slope by
Lemma \ref{prelagrangiano}.
 
Suppose that $T \in {\mathcal T}_{\gamma}$ contains a convex torus $T$ with finite 
slope. First, we show that  $T \in {\mathcal T}_{\gamma}$ also contains a convex 
torus $T'$ with infinite slope. In fact, by hypothesis there is a vertical 
Legendrian curve $L \subset M$ with twisting number $0$, hence we obtain a
convex torus with infinite slope by first isotoping $T$ so that it becomes a 
convex torus $T'$ with vertical ruling disjoint from $L$. Then, if $T'$ has not
infinite slope already, by attaching the 
bypasses along $T'$ coming from a convex annulus between $L$ and a Legendrian 
ruling curve of $T'$. This operation produces a convex torus $T''$ parallel to 
$T$ with infinite slope. Once we have a convex torus with finite slope $T$ and a 
convex torus with infinite slope $T''$, we can suppose by isotopy discretisation 
that they are disjoint, so they bound a tight thickened torus with different 
boundary slopes. By \cite{honda:1}, Corollary 4.8, such thickened torus contains 
a pre-Lagrangian torus.
\end{proof}

\begin{proof}[Proof of Proposition \ref{conto}]
Take a curve $\gamma' \subset T^2 \setminus D= \Sigma$ isotopic to $\gamma$  which realises the 
minimum of the intersection with the dividing set $\Gamma$. We can identify
$\gamma'$ with its image under the section $\sigma$ and make it 
Legendrian. Since $(M \setminus V, \, \xi|_{M \setminus V}) \cong (\Sigma \times S^1, \, \xi_{\Gamma_{\Sigma}})$ is 
$S^1$--invariant, $T_0= \pi^{-1}(\gamma')$ is a standard torus with division number
$\frac{1}{2} | \gamma \cap \Gamma |$ if $\gamma' \cap \Gamma \neq \emptyset$, 
or a pre-Lagrangian torus if $\gamma' \cap \Gamma = \emptyset$, therefore 
\[\min \limits_{T \in {\mathcal T}_{\gamma}} \hbox{div} (T) \leq \frac{1}{2} | \gamma \cap \Gamma |.\] 
Suppose by contradiction that there exists a convex torus 
$T_1 \in {\mathcal T}_{\gamma}$  with either $\hbox{div} (T_1) < \frac{1}{2} |\gamma \cap 
\Gamma|$ and $|\gamma \cap \Gamma|>2$, or with slope different from infinity and $|\gamma \cap \Gamma|=2$. 
By Isotopy Discretisation \cite{gluing}, Lemma 3.10, we can find 
a finite family of disjoint convex tori
$T_0=T^{(0)}, \ldots ,T^{(n)}=T_1$ such that, for any $i=0, \ldots ,n-1$, 
 $T^{(i+1)}$ is obtained from $T^{(i)}$ by the attachment of a single bypass. 
In particular, they bound a layer $N_i$ diffeomorphic to $T^2 \times I$.  
If $T^{(n)}=T_1$ has finite slope, we can assume that it is the first torus in
the family with that property.

For any $i$ such that $T^{(i)}$ has infinite slope there is a Seifert 
fibration $\pi_i\co M \setminus T^{(i)} 
\to S^1 \times I$ with one singular fibre, a neighbourhood $V_i \subset M 
\setminus N_i$ of the singular fibre such
that $- \partial (M \setminus V_i)$ has infinite slope, and a collar 
$C_i= \pi_i (N_i)$ of a boundary component of $\Sigma_0$ such that
\begin{align} 
\nonumber & \pi_i |_{N_i}\co N_i \to C_i \\
\nonumber & \pi_i |_{M \setminus (N_i \cup V_i)}\co M \setminus (N_i \cup V_i) \to \Sigma_0 \setminus C_i
\end{align} 
are $S^1$--bundles. 
We choose sections $\sigma_i\co \Sigma_0 \to M \setminus (T^{(i)} \cup V_i)$ 
so that: 
\begin{enumerate}
\item $\sigma_i(\Sigma_0)$
is a convex $\# \Gamma$--minimising surface with Legendrian boundary
 denoted by $\Sigma_{(i)}$. 
\item $\Sigma_{(i)} \cap T^{(i+1)}$ is a Legendrian curve.
\item $\sigma_i$ extends to a section $\overline{\sigma}_i\co \Sigma \to M \setminus V_i$.
\item $\overline{\sigma}_0 = \sigma$.
\end{enumerate}

Define $\overline{\Sigma}_{(i)}= \overline{\sigma}_i(\Sigma)$ and identify 
$\Gamma_{\overline{\Sigma}_{(i)}}$ with a multicurve on $\Sigma$ using $\pi_i$. 

We claim that, for any $i$, $\Gamma_{\overline{\Sigma}_{(i)}}$  differs from 
$\Gamma_{\Sigma}$ by a number of curves isotopic to $\gamma$ or by Dehn twists around $\gamma$. 
The proof is by induction on $i$.   If $i=0$ the claim is true because 
$\Sigma = \overline{\Sigma}_{(0)}$. Now suppose the claim true 
for a fixed $i$. 
Let $\sigma_{i+1}'\co \Sigma_0 \to M \setminus (V_i \cup T^{(i+1)})$ be the section which extends to the section
$\overline{\sigma}_i\co \Sigma \to M \setminus V_i$. We denote $\sigma_{i+1}'(\Sigma_0)$ by $\Sigma_{i+1}'$ and 
$\overline{\sigma}_{i+1}'(\Sigma)$ by $\overline{\Sigma}_{i+1}'$. By properties (1) and (2) of $\sigma_i$, 
$\Sigma_{i+1}'$ is a convex $\# \Gamma$--minimising surface with Legendrian boundary, then by 
\cite{honda:2}, Lemma 4.1  the $S^1$--invariant contact manifolds
$(M \setminus (T_{(i+1)} \cup V_i), \, \xi |_{M \setminus (T_{(i+1)} \cup V_i)})$ and
$(\Sigma_0 \times S^1, \, \xi_{\Gamma_{\Sigma_{i+1}'}})$ are contactomorphic. Analogously, 
$(M \setminus (T_{(i+1)} \cup V_{i+1}), \, \xi |_{M \setminus (T_{(i+1)} \cup V_{i+1})})$ is contactomorphic to 
$(\Sigma_0 \times S^1, \, \xi_{\Gamma_{\Sigma_{i+1}}})$.
These contactomorphisms give presentations of $(M \setminus T_{(i+1)}, \, \xi |_{M \setminus T_{(i+1)}})$
as negative contact surgery on $(T^2 \times I, \, \xi_{\Gamma_{\Sigma_{i+1}'}})$ and 
$(T^2 \times I, \, \xi_{\Gamma_{\Sigma_{i+1}}})$ respectively, therefore  
by Lemma \ref{confronto}, 
 $\widehat{\Gamma}_{\Sigma_{i+1}}$ is diffeomorphic to $\widehat{\Gamma}_{\Sigma_{i+1}'}$. By construction,
$\Gamma_{\Sigma_i}$ and $\Gamma_{\Sigma_{i+1}'}$ extend to the same 
multicurve $\Gamma_{\overline{\Sigma}_i}$ on $\overline{\Sigma}$, so the claim is proved. 

Suppose now that $T_1$ has infinite slope and $\hbox{div}(T_1)< \hbox{div}(T_0)$:
then the geometric intersection $|\gamma \cap  \Gamma_{\overline{\Sigma}_{(n)}}|$ is lesser than
the the geometric intersection $|\gamma \cap \Gamma_{\Sigma_{(0)}}|=|\gamma \cap \Gamma|$. This is a
contradiction because, by the claim, $\Gamma_{\overline{\Sigma}_{(0)}}$ and 
$\Gamma_{\overline{\Sigma}_{(n)}}$ differ only by Dehn twists along $\gamma$ or by the number of 
curves isotopic to $\gamma$. 

If the slope of $T_1$ is not infinity, attaching the bypasses coming from a 
vertical annulus $A \subset M \setminus N_{n-1}$ between $T_{(n)}$ and $T_{(n-1)}$  we find a layer 
$N_{n} \cong T^2 \times I$ so that $N= N_{n-1} \cup 
N_n$ has minimal boundary and infinite boundary slopes. $N$ is rotative because
it has infinite boundary slopes, but $T_{(n)} \subset N$ has finite slope, so  
by \cite{honda:1}, Lemma 5.7,  the dividing set of a $\# \Gamma$--minimising
 section of $N$ contains no arcs with endpoints on different boundary
components.  We
can complete the section in $N$ to a section $\Sigma_{(n-1)}'= \sigma_{n-1}'(\Sigma_0)$ which has 
no dividing arcs with endpoints on different boundary components.
By Lemma \ref{confronto}, $ \widehat{\Gamma}_{\Sigma_{(n-1)}'}$ is isomorphic to 
$\widehat{\Gamma}_{\Sigma_{(n-1)}}$
and by the claim $ \widehat{\Gamma}_{\Sigma_{(n-1)}}$ glues to the same dividing set on $T^2$ as
 $ \widehat{\Gamma}_{\Sigma_{(0)}}$. This is a contradiction because $ \widehat{\Gamma}_{\Sigma_{(n-1)}'}$ 
has a curve isotopic to $\gamma$ and $ \widehat{\Gamma}_{\Sigma_{(0)}}$ does not.
\end{proof}

\begin{theorem}\label{paperino}
Let $\Gamma_{\Sigma}$ and $\Gamma_{\Sigma}'$ be two tight abstract dividing sets on the punctured torus
$\Sigma$ such that $\# \Gamma_{\Sigma} \cap \partial \Sigma = \# \Gamma_{\Sigma}' \cap \partial \Sigma =2$ and without boundary parallel dividing 
arcs. Denote by $\Gamma$ and $\Gamma'$ their
completion.
The tight contact structures $\xi_{\Gamma_{\Sigma}}(\eta)$ and $\xi_{\Gamma_{\Sigma}'}(\eta')$ on   $M(e_0, r)$ 
are isotopic if and only if $\Gamma$ is isotopic to $\Gamma'$, and $\eta$ is 
isotopic to $\eta'$.
\end{theorem}
\begin{proof}
If $\Gamma$ is isotopic to $\Gamma'$ and $\eta$ is isotopic to $\eta'$, then 
$\xi_{\Gamma_{\Sigma}}(\eta)$ is isotopic to $\xi_{\Gamma_{\Sigma}'}(\eta')$
by Proposition \ref{ultima}.

Let now $\xi_{\Gamma_{\Sigma}}(\eta)$ and $\xi_{\Gamma_{\Sigma}'}(\eta')$
be isotopic tight contact structures.
By Proposition \ref{conto} for any simple closed curve $\gamma \subset T^2$ we have 
$| \gamma \cap \Gamma | = | \gamma \cap \Gamma' |$, therefore $\Gamma_{\Sigma}$ is isotopic to $\Gamma_{\Sigma}'$. Suppose now 
that there are tight contact structures $\eta_0$ and $\eta_1$ on $D^2 \times S^1$ such that 
$\xi_{\Gamma_{\Sigma}}(\eta)$ is isotopic to $\xi_{\Gamma_{\Sigma}}(\eta')$: then there exist isotopic convex 
vertical tori with infinite slope $T_0$, $T_1 \subset M$, and, for $i=0,1$, Seifert
fibrations $\pi_i\co M \setminus T_i \to S^1 \times I$, and neighbourhoods $V_i$ of the 
singular fibre such that $\xi|_{V_i} \cong \eta_i$. 

By isotopy discretisation \cite{gluing}, Lemma 3.10, there is a finite 
sequence of convex tori with infinite slope $T_0=T^{(0)}, \ldots ,T^{(n)}=T_1$ such that, 
for $i=0, \ldots ,n-1$, $T^{(i)}$ and $T^{(i+1)}$ bound  $N_i$ diffeomorphic to $T^2 \times I$. 
For any $i$, we can modify the
Seifert fibration on $M$ so that the singular fibre is contained in $M \setminus N_i$,
and find a neighbourhood of the singular fibre $V_i'$ contained in $M \setminus N_i$.
By Lemma \ref{confronto} applied to $M \setminus T^{(0)}$, $\eta_0 =\xi|_{V_0}$ is isotopic to 
$\xi|_{V_0'}$, by  Lemma \ref{confronto} applied to $M \setminus T^{(i)}$, $\xi|_{V_i'}$ is 
isotopic to $\xi|_{V_{i+1}'}$, and by  Lemma \ref{confronto} applied to 
$M \setminus T^{(n)}$,
 $\xi|_{V_n'}$ is isotopic to $\xi|_{V_1} \cong \eta_1$.
\end{proof}

%%%%%%%%%%%%%%%%%%%%%%%%%%%%%%%%%%%%%%%%%%%%%%%%%%%%%%%%%%%%%%%%%%%%%%%%%%%%%%
%%%%%%%%%%%%%%%%%%%%%%%%%%%%%%%%%%%%%%%%%%%%%%%%%%%%%%%%%%%%%%%%%%%%%%%%%%%%%%
\section{Exceptional tight contact structures}
In this section we prove tightness for the candidate tight contact structures 
with $\# \Gamma =1$. The proof of tightness for this class of contact structures uses 
a purely topological and three dimensional technique known as 
{\em state traversal}, introduced by Honda in \cite{honda:2}. 

\subsection{State traversal}
Let $(M, \, \xi)$ be a contact manifold and $W \subset M$ be a properly embedded 
incompressible surface. We will assume $W$ is convex. The 
contact manifold $(M \setminus W, \, \xi|_{M \setminus W})$ will be called a {\em state}, and the 
surface  $W$ a {\em wall}. In general both $W$ and $M \setminus W$  
could be 
disconnected. A state is said tight if $\xi|_{M \setminus W}$ is tight. The boundary of 
$M \setminus W$ consists of two copies of $W$: $W_+$ and $W_-$. A {\em state transition} 
consists
of detaching a collar of $W_-$ and attaching it to $W_+$, or vice versa, so 
that
$\Gamma_W$ is changed by a bypass attachment.
We observe that a state transition corresponds to moving $W$ inside $M$ by an
isotopy.
\begin{theorem}[\cite{honda:2}, section 2.3.1 or \cite{gluing}, Theorem 3.5]
If the initial state $(M \setminus W, \, \xi|_{M \setminus W})$ is tight, and all the states reached
 from it in a finite number of state transitions are also tight, then
the contact manifold $(M, \, \xi)$ is tight.
\end{theorem}

The set of states that can be reached from the initial state in a finite number
of bypass attachments is a complete isotopy invariant of $\xi$ in the following
sense.
\begin{theorem}[Corollary of \cite{gluing}, Theorem 3.1] 
Let $\xi_1$ and $\xi_2$ be two tight contact structures on $M$, and let $W \subset M$ be a 
properly embedded incompressible convex surface. Let
${\mathcal C}(\xi_i)$ be the set of isotopy classes relative to the boundary 
of all the states reached from the initial state  $(M \setminus W, \, \xi|_{M \setminus W})$
in a finite number of state transitions. Then $\xi_1$ is isotopic to 
$\xi_2$ if and only if ${\mathcal C}(\xi_1) = {\mathcal C}(\xi_2)$.
\end{theorem} 

Let $\Gamma_{\Sigma}$ be an abstract dividing set on the one--punctured torus $\Sigma$ with
$\# \Gamma_{\Sigma}=1$. In order to apply the state traversal to $(M(e_0,r), \, \xi_{\Gamma_{\Sigma}}(\eta))$, 
we give an alternative description of this contact structure.
By $\xi_l$, with $l=2$ or $l=-2$, we denote the tight contact structure with 
infinite boundary slopes and twisting $\pi$ on $T^2 \times I$ such that $\xi_l$ has 
relative Euler class $e(\xi_l)=  \left ( \begin{array}{c} 0 \\ l \end{array} 
\right )$, ($\xi_1^{\pm}$ in the notations
of \cite{honda:1}, Lemma 5.2) and by $(M', \, \xi_l(\eta))$ we denote the contact 
manifold obtained by contact $(- \frac{1}{r})$--surgery
along a vertical Legendrian curve in $(T^2 \times I, \, \xi_l)$ with  twisting number 
$0$.
\begin{figure}[ht!]
\centering
\psfrag{a}{\footnotesize $T^2 \times \{ 0 \}$}
\psfrag{b}{\footnotesize $T^2 \times \{ 1 \}$}
\includegraphics[width=6cm]{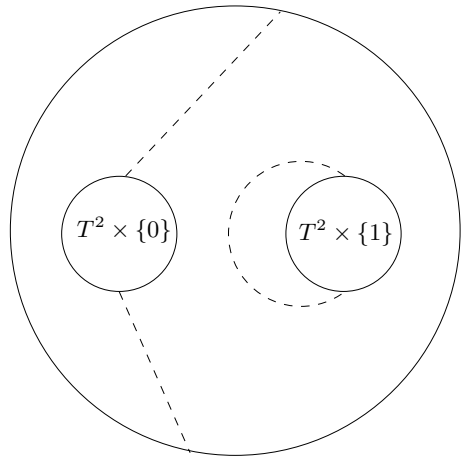}
\caption{The dividing set $\Gamma_{\Sigma_0}$}
\label{ultimafigura.fig} 
\end{figure}
\begin{lemma}
After gluing $T_1$ to $T_0$ with the map $\left ( \begin{matrix}
1 & 0 \\
-e_0 & 1 
\end{matrix} \right )$, we obtain the contact manifold 
$(M(e_0,r), \, \xi_{\Gamma_{\Sigma}}(\eta))$ with $\# \Gamma_{\Sigma}=1$ and $\langle e(\xi_{\Gamma_{\Sigma}}), \, \Sigma \rangle = l$. 
\end{lemma}
\begin{proof}
To prove that we obtain a contact structure isotopic to $\xi_{\Gamma_{\Sigma}}(\eta)$ it is 
enough to show that we obtain a contact manifold whose background is
isotopic to the background of $(M(e_0,r), \, \xi_{\Gamma_{\Sigma}}(\eta))$. We prove the isotopy 
between the backgrounds by showing that they induce isotopic dividing sets
on convex $\# \Gamma$--minimising sections of $\Sigma \times S^1$, see \cite{honda:2} Lemma 4.1.
 
The background of $(M', \, \xi_l(\eta))$ is the contact manifold 
$(\Sigma_0 \times S^1, \, \xi_{\Gamma_{\Sigma_0}})$ where $\Sigma_0$ is a pair of pants and $\xi_{\Gamma_{\Sigma_0}}$ is the 
tight contact structure which is $S^1$--invariant over the abstract dividing 
set described in Figure \ref{ultimafigura.fig}. $\Gamma_{\Sigma_0}$ consists of two arcs 
joining two boundary components of $\Sigma_0$ and a boundary parallel arc with both
end-points on the third component of $\partial \Sigma_0$. 

Consider a $\# \Gamma$--minimising section $\Sigma_0'$ of $\Sigma_0 \times S^1$ so that, after gluing 
two boundary components of $\Sigma_0 \times S^1$ as prescribed by the statement, we obtain 
a section $\Sigma'$ of $\Sigma \times S^1$. By \cite{honda:2} Lemma 4.5 $\Gamma_{\Sigma_0'}$ is isotopic 
to $\Gamma_{\Sigma_0}$, therefore $\Gamma_{\Sigma'}$ is isotopic to $\Gamma_{\Sigma}$. 

\end{proof}
For the rest of the 
section we fix the notation $M=M(e_0,r)$ and $\xi = \xi_l(\eta) = \xi_{\Gamma_{\Sigma}}(\eta)$.

\begin{lemma}\label{infinito}
Let $W_0, W_1 \subset (M, \, \xi)$ be  vertical incompressible disjoint convex tori with 
finite slopes and let $N \subset M$ be the thickened torus bounded by $W_0$ and
$W_1$. If $\partial (M \setminus N)= W_1 - W_0$, the slope of $W_1$ is $s_1= \frac{p}{q}$, and the 
slope of $W_0$ is $s_0 \leq \frac{p}{q} -e_0$, then there is a vertical Legendrian 
curve $L \subset M \setminus N$  with twisting number $tb(L)=0$.
\end{lemma}
\begin{proof}
Take a properly embedded convex vertical annulus with Legendrian boundary
$A \subset (M \setminus N)$ and attach all the possible bypasses it carries to $W_0$ and 
$W_1$. If in the process we get a torus with infinite slope
we are done, otherwise we end with two convex tori $W_0'$ and $W_1'$ parallel 
to $W_0$ and $W_1$ respectively with slopes $s_0' \leq s_0 \leq \frac{p}{q} -e_0$ 
and $s_1' \geq \frac{p}{q}$, and such that the vertical annulus $A'$ between $W_i'$ 
and $W_{i+1}'$ carries no bypasses. Let $N'$ be the layer diffeomorphic to
$T^2 \times I$ between $W_0'$ and $W_1'$. 
If we cut $M \setminus N'$ along $A'$ and round the edges, we get 
a solid torus with boundary slope $s \geq s_1' -s_0' -1 \geq 0$. On the other hand, if
we make the singular fibre $F$ Legendrian with very low twisting number $n$, 
and remove a standard neighbourhood $\nu F$, we get
slope $\frac{- n \beta - \beta'}{n \alpha + \alpha'}$ on $- \partial (M' \setminus \nu F)$. 
Taking the limit for $n$ going to infinity, we see that this 
slope is negative for $n$ small enough, therefore, by \cite{honda:1}, 
Proposition 4.16, there 
is an intermediate torus  in $M \setminus (N' \cup A' \cup \nu F)$   with infinite slope.
\end{proof}

If $W \subset  (M, \, \xi)$ is a vertical incompressible convex torus with finite 
slope, applying Lemma \ref{infinito} with $W_0=W_1=W$, we find a vertical 
Legendrian curve $L$ with twisting number zero in $M \setminus W$. Attaching the 
bypasses coming from $L$ to $W$ on either sides,  we can  engulf $W$ in a 
rotative $T^2 \times [-1,1]$ with infinite boundary slopes such that $W=T^2 \times \{ 0 \}$. 
By \cite{corrige}, Theorem 2.2, we need to consider only 
transitions between states with  minimal boundary, provided that the walls can 
be engulfed into rotative thickened tori. In the present situation, we have to
consider transitions between states with finite boundary slope and minimal 
boundary, or between states with infinite boundary slopes. 

\subsection{Analysis of the states}
Before performing the state traversal we 
analyse the possible states. We observe that the Seifert fibration on $M$ can
be isotoped so that $W$ becomes a fibred torus. Consequently there is an 
induced Seifert fibration on $M \setminus W \cong M'$. 

Let $\zeta_l'$ be the  minimally 
twisting tight contact structure on $T^2 \times [0, \frac 12 ]$ with boundary slopes 
$s_0 = \frac{p}{q} -e_0$, $s_{\frac 12 }= \infty$, minimal boundary and relative
 Euler class $\pm \left ( \begin{array}{c} -q \\ -1-p+e_0q \end{array} \right )$. 
Let $\zeta_l''$ be the minimally 
twisting tight contact structure  on $T^2 \times [ \frac 12 ,1]$ with boundary slopes
$s_{\frac{1}{2}} = \infty$, $s_1= \frac{p}{q}$, minimal boundary and relative
 Euler class $\pm \left ( \begin{array}{c} q \\ 1+p \end{array} \right )$. Here 
the signs of the relative Euler 
classes are chosen accordingly to the sign of $l$. All the basic slices in the 
decomposition of $\zeta_l'$ and $\zeta_l''$ have the same sign.
We denote by $(M'',\, \zeta_l'(\eta))$ the contact manifold 
constructed by contact $(- \frac{1}{r})$-surgery on a vertical 
Legendrian curve $L$ with twisting number $0$ in 
$(T^2 \times [0, \frac 12 ], \, \zeta_l')$, and by $(M', \, \xi_l'(\eta))$ the contact manifold  
obtained by gluing $(T^2 \times [ \frac 12 ,1], \, \zeta_{-l}'')$ to $(M'', \, \zeta_l'(\eta))$.

To the contact manifold $(D^2 \times S^1, \eta)$ we associate the set of 
numerical invariants $(r_0, \ldots ,r_k)$ defined as 
\[r_i= \# \{ \hbox{positive basic slices in } N_i \} - \# 
\{ \hbox{negative basic slices in } N_i \}\]
where $N_i$ is the $(i+1)$--th continued fraction block in the basic slices 
decomposition of $\eta$.
The classification of  tight contact structures on 
solid tori \cite[Theorem 2.3]{honda:1}
 implies the following proposition.
\begin{prop}\label{comoda} 
For any slope $s$, let $\Gamma_s$ denote the tight abstract dividing set on $T^2$ 
with 
slope $s$ and $\# \Gamma_s=2$. Then the number of continued fraction blocks  and 
the slopes of the borders between continued fraction blocks of any tight 
contact structure $\eta \in \hbox{Tight}(D^2 \times S^1, \Gamma_s)$ depend only on $s$ and are 
independent of $\eta$. Moreover, the map $\hbox{Tight}(D^2 \times S^1, \Gamma_s) \to \Z^{n+1}$ 
given by $\eta \mapsto (r_0, \ldots ,r_k)$ is injective.
\end{prop} 
If $\Gamma_s=A(r)^{-1} \Gamma_{\infty}$, then the number of continued fraction 
blocks of $\eta \in \hbox{Tight}(D^2 \times S^1, \Gamma_s)$ is equal to the 
number of the coefficient in the continued fraction expansion 
$- \frac{1}{r} =[d_0, \ldots ,d_n]$, moreover $|r_0| \leq |d_0|$, and $|r_i| \leq |d_i|-1$
for $i>0$. To classify the contact manifolds $(M', \, \xi_l'(\eta))$ we need the 
following lemmas.

\begin{lemma} \label{pre-modello} 
Let $(T^2 \times [0, \frac 12 ],  \, \zeta_l')$ be a basic slice with boundary slopes $-n$ 
and $\infty$, and let $(T^2 \times [\frac 12 ,1], \, \zeta_{-l}'')$ be a basic slice with 
boundary slopes $\infty$ and $1$ with opposite sign. We call 
$(T^2 \times [0,1] \setminus V, \, \xi_l')$ the contact manifold obtained from 
$(T^2 \times [0, \frac 12 ],  \, \zeta_l') \cup_{T^2 \times \{ \frac 12 \}} 
(T^2 \times [\frac 12 ,1], \, \zeta_{-l}'')$ by removing a standard 
neighbourhood $V$ of a vertical Legendrian divide of $T^2 \times \{ \frac 12 \}$. Then
$(T^2 \times [0,1] \setminus V, \, \xi_l')$ is isomorphic to $(T^2 \times [0, \frac{1}{2}] \setminus U, \,
\zeta_l'|_{T^2 \times [0, \frac{1}{2}] \setminus U})$, where $U$ is a standard neighbourhood 
of a vertical Legendrian ruling curve of a standard torus parallel to
$T^2 \times \{ 0 \}$ and contained in its invariant neighbourhood.
\end{lemma}
\begin{proof}
Up to isotopy we can assume that $T^2 \times [\frac 12 - \epsilon, \frac 12 + \epsilon]$ is an
invariant neighbourhood of $T^2 \times \{ \frac 12 \}$ in $(T^2 \times [0, \frac 12 ],  \, 
\zeta_l') \cup  (T^2 \times [\frac 12 ,1], \, \zeta_{-l}'')$, and that $V$ is contained in it. 
Clearly $\xi_l'|_{T^2 \times [0, \frac 12 - \epsilon]}$ is isomorphic to $\zeta_l'$ and 
$\xi_l'|_{T^2 \times [\frac 12 + \epsilon, 1]}$ is isomorphic to $\zeta_{-l}''$. By \cite[Lemma 4.1]{honda:2}  $\xi_l'|_{T^2 \times [\frac 12 - \epsilon, \frac 12 + \epsilon] \setminus V}$ is $S^1$--invariant, and 
the dividing set of a convex $\# \Gamma$--minimising horizontal section $\Sigma$ with 
Legendrian boundary is as in Figure \ref{divisione.fig}.

\begin{figure}[ht!] 
\centering
\includegraphics[width=4cm]{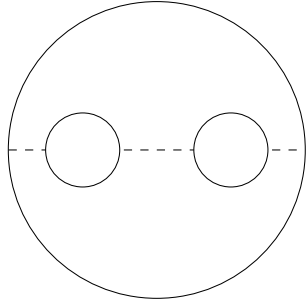}
\caption{The dividing set on $\Sigma$ and $\Sigma'$}
\label{divisione.fig}
\end{figure}

We decompose $(T^2 \times [0, \frac{1}{2}] \setminus U, \, \zeta_l'|_{T^2 \times [0, \frac{1}{2}] \setminus U})$
into three pieces which are isomorphic to $(T^2 \times [0, \frac 12 - \epsilon], \,
 \xi_l'|_{T^2 \times [0, \frac 12 - \epsilon]})$, $(T^2 \times [\frac 12 + \epsilon, 1], \, 
\xi_l'|_{T^2 \times [\frac 12 + \epsilon, 1]})$, and $(T^2 \times [\frac 12 - \epsilon, \frac 12 + \epsilon] \setminus V, \,
\xi_l'|_{T^2 \times [\frac 12 - \epsilon, \frac 12 + \epsilon] \setminus V})$ respectively. Thicken 
$U \in (T^2 \times [0, \frac{1}{2}], \, \zeta_l')$ to $U'$ by attaching the bypasses
coming from a vertical annulus between a Legendrian ruling curve of $\partial U$ and
a  Legendrian divide of $T^2 \times \{ \frac 12 \}$ so that $\partial U'$ has infinite slope.
In a similar way find a collar $C$ of $T^2 \times \{ 0 \}$ in 
$(T^2 \times [0, \frac 12 ] \setminus U', \, \zeta_l'|_{T^2 \times [0, \frac 12 ] \setminus U'})$ so that $C$ has
boundary slopes $-n$ and $\infty$. The isomorphism between $(T^2 \times [0,1] \setminus V, \,
\xi_l')$ and $(T^2 \times [0, \frac 12 ] \setminus U, \, \zeta_l'|_{T^2 \times [0, \frac 12 ] \setminus U})$
identifies $T^2 \times [0, \frac 12 - \epsilon]$ to $C$, $T^2 \times [\frac 12 - \epsilon, \frac 12 + \epsilon] 
\setminus V$ to $T^2 \times [0, \frac 12 ] \setminus (U' \cup C)$, and $T^2 \times [\frac 12 + \epsilon, 1]$ to 
$U' \setminus U$. 

$(T^2 \times [0, \frac 12 ] \setminus (U' \cup C), \, \zeta_l'|_{T^2 \times [0, \frac 12 ] \setminus (U' \cup C)})$ has
infinite boundary slopes, therefore by \cite[Lemma 4.1]{honda:2} it is
 $S^1$--invariant. Let $\Sigma'$ be a convex $\# \Gamma$--minimising horizontal section 
 with Legendrian boundary. The dividing sets of $\Sigma'$ cannot contain any 
boundary parallel dividing arc, otherwise such an arc would produce a bypass 
attached horizontally to  $T^2 \times [0, \frac 12 ] \setminus (U' \cup C)$. The attachment of 
this bypass would give a convex torus with slope zero in 
$(T^2 \times [0, \frac 12 ], \, \zeta_l')$, contradicting either tightness or minimal 
twisting. Thus the dividing sets of $\Sigma'$ is forced to be as in Figure 
\ref{divisione.fig}, therefore $(T^2 \times [0, \frac 12 ] \setminus (U' \cup C), \, \zeta_l'|_{T^2 \times 
[0, \frac 12 ] \setminus (U' \cup C)})$ is isomorphic to 
$(T^2 \times [\frac 12 - \epsilon, \frac 12 + \epsilon] \setminus V, \, 
\xi_l'|_{T^2 \times [\frac 12 - \epsilon, \frac 12 + \epsilon] \setminus V})$ by \cite[Lemma 4.1]{honda:2}
because they induce diffeomorphic dividing sets on  the convex 
$\# \Gamma$--minimising horizontal sections $\Sigma$ and $\Sigma'$.

$\zeta_l'|_{C}$ is isomorphic to $\zeta_l' \cong \xi_l'|_{T^2 \times [0, \frac 12 - \epsilon ]}$ because their 
relative Euler classes have the same evaluation on a vertical annulus.  In a 
similar way, the relative Euler class of $\zeta_l'|_{U' \setminus U}$ evaluates on a vertical 
annulus as the relative Euler class of $\zeta_l'$, therefore it evaluates as the 
opposite of the relative Euler class of 
$\zeta_{-l}'' \cong \xi_l'|_{T^2 \times [\frac 12 + \epsilon , 1]}$. The change of sign is due to the fact
that, in evaluating the relative Euler class of $\zeta_l'|_{U' \setminus U}$, the boundary 
component $\partial U$ is oriented as $T^2 \times \{ 0 \}$, i.~e. by the inward normal.
On the other hand, the isomorphism maps $\partial U$ to $T^2 \times \{ 1 \}$, which is 
oriented by the outward normal. This change of orientation on the boundary
forces the orientation of the vertical annulus to change too, in order to
keep the global orientation unchanged. 
\end{proof}
The contact manifold $(T^2 \times [0, \frac 12 ],  \, \zeta_l') \cup (T^2 \times [\frac 12 ,1], \, 
\zeta_{-l}'')$ is overtwisted because it does not satisfy the conditions of the 
Gluing Theorem \cite{honda:1}, Theorem 4.25. However, it become tight when we 
remove $V$, as Lemma \ref{pre-modello} shows.

\begin{lemma} \label{modello}
Let $(T^2 \times [0,1] \setminus V , \, \xi_l')$ be a contact manifold as in Lemma 
\ref{pre-modello}. 
Then  there is a convex annulus $A \subset T^2 \times I \setminus V$ whose boundary
consists of vertical Legendrian ruling curves  of $T^2 \times \{ 0 \}$ 
and $T^2 \times \{ 1 \}$ such that its dividing set $\Gamma_A$ has no boundary parallel 
dividing curves, and $\xi_l'|_{T^2 \times I \setminus (V \cup A)}$  is isotopic to $\zeta_l'$.
\end{lemma}

\begin{proof}
The annulus $A$ can be easily found in the contact manifold
$(T^2 \times [0, \frac 12],\, \zeta_l')$, which is isomorphic to $(T^2 \times [0,1] \setminus V, \, \xi_l')$
by Lemma \ref{pre-modello}. It is a convex vertical annulus with Legendrian 
boundary between $T^2 \times \{ 0 \}$ and $\partial U$ contained in an invariant collar of 
$T^2 \times \{ 0 \}$. The invariance of the collar implies that the annulus can carry
no bypass.
\end{proof}

\begin{lemma} \label{stato1}
If $e_0 \geq 2$, then $(M', \, \xi_l'(\eta))$ is tight for any $\eta$ and any $l$. If 
$e_0=1$, then $(M', \, \xi_l'(\eta))$ is tight if and only if $r_0= \frac{ld_0}{2}$.
Moreover, $\xi_{l'}'(\eta')$ is isotopic to $\xi_l'(\eta)$ if and only if either $l=l'$ and 
$\eta$ is isotopic to $\eta'$ or $s_0 \in \Z$ and 
\begin{enumerate}
\item  $l'=-l$ and $r_0'= r_0+l$ when $e_0=2$,
\item $l'=-l$, $r_0=-r_0'$ and $r_1'= r_1+l$ when $e_0=1$.
\end{enumerate}
\end{lemma}
\begin{proof}
After acting, if necessary, on $M'$ by a self-diffeomorphism
supported outside the surgery and preserving the Seifert fibration, 
we can assume $s_1 = \frac{p}{q} \in (0,1]$.
There is a unique universally tight contact structure on $D^2 \times S^1$ with 
boundary slope $-s_1$ which can be glued to $T^2 \times [\frac{1}{2} ,1]$
along $- T^2 \times \{ 1 \}$ with the identity map to give the tight contact structure
with infinite boundary slope on $T^2 \times [\frac 12 ,1] \cup D^2 \times S^1 \cong D^2 \times S^1$.
After filling $T^2 \times \{ 1 \}$ in $(M', \, \xi_l'(\eta))$ in this way, we obtain a 
contact structure on $D^2 \times S^1$ still denoted by $\xi_l'(\eta)$ which can be 
decomposed
 as $(D^2 \times S^1, \, \xi_l'(\eta)) \cong (T^2 \times [0, \frac{1}{2}], \, \zeta_l'') \cup_{A(r)} (D^2 \times S^1, \eta)$.

We start studying tightness for this bigger contact manifold. Since 
both $(T^2 \times [0, \frac{1}{2}], \, \zeta_l'')$ and $(D^2 \times S^1, \eta)$ are tight, the 
Gluing Theorem \cite{honda:1}, Theorem 4.25,   
gives necessary and sufficient conditions for $\xi_l(\eta)$ to
be tight. The application of the Gluing Theorem for thickened tori to 
solid tori is possible
because \cite{honda:1}, Propositions 4.15, 4.17, and 4.18 give an 
identification between isotopy classes of tight contact structures on 
$D^2 \times S^1$ and isotopy
classes of tight minimally twisting contact structures on $T^2 \times I$.

Let $\frac{p}{q} -e_0=s_0 > s_{\frac{1}{2n}} > \ldots > s_{\frac{n-1}{2n}}> s_{\frac{1}{2}}
= \infty$ be the sequence of the boundary slopes of a minimal basic slices 
decomposition of $(T^2 \times [0, \frac{1}{2}], \, \zeta_l')$,
and let $\infty =s_{\frac{1}{2}}' > s_{\frac{m-1}{2m}}' > \ldots > s_{\frac{1}{2m}}' > s_0'$ be
the sequence of the boundary slopes of a minimal basic slices decomposition 
of $(D^2 \times S^1, \eta)$ computed with respect to the basis of 
$T^2 \times [0, \frac{1}{2}]$.
 The contact 
structure $\xi_l(\eta)$ is tight if and only if either $s_0 > \ldots > s_{\frac{n-1}{2n}} >
\infty > s_{\frac{m-1}{2m}}' > \ldots > s_0'$ is the shortest sequence of slopes between $s_0$ and
 $s_0'$, or there exist $s_{\frac{k}{2n}} < \infty < s_{\frac{k'}{2m}}'$ joined by an edge in 
the Farey tessellation 
of ${\mathbb H}^2$ such that the basic slices between them  have all the same 
signs.  

If $\frac{p}{q} -e_0 <0$, then the shortest sequence of slopes between 
$s_0$ and $s_0'$ needs to go
through $\infty$, because all the $s_{\frac{i}{2n}}$ are negative and all the
$s_{\frac{j}{2m}}'$ are positive, therefore $(D^2 \times S^1, \, \xi_l'(\eta))$ is tight for any $l$ 
and any $\eta$.
If $\frac{p}{q} -e_0 =0$, which means $\frac{p}{q} =1$, and $e_0=1$, there is an
edge in the Farey Tessellation joining $0$ with $- \frac{1}{d_0+1}$, so
$\infty$ is not a border between basic slices. In this case 
$(D^2 \times S^1, \, \xi_l'(\eta))$  is tight if and only if  $(T^2 \times [0, \frac 12 ], \, \zeta_l'')$ 
with boundary slopes $s_0=0$
and $s_1= \infty$ glues with the first continued fraction block of $\eta$ with boundary 
slopes $s_{\frac 12}'= \infty$ and $s_{\frac{m -d_0-1}{2m}}'= - \frac{1}{d_0+1}$ to give 
a basic slice. This 
happens if and only if $l=2$ and $r_0=d_0$, or $l=-2$ and $r_0=-d_0$. When 
$(D^2 \times S^1, \, \xi_l'(\eta))$ is not tight, then $(M', \, \xi_l'(\eta))$ is not tight either, 
because Lemma \ref{modello} gives a contact embedding of $(D^2 \times S^1, \, \xi_l'(\eta))$ 
into $(M', \, \xi_l'(\eta))$.

In order to determine whether $(M', \, \xi_l'(\eta))$ is isotopic to 
$(M', \, \xi_{l'}'(\eta'))$, we again study the problem in $D^2 \times S^1$ first. In fact, if 
$(M', \, \xi_l'(\eta))$ is isotopic to $(M', \, \xi_{l'}'(\eta'))$, then $(D^2 \times S^1, \, \xi_l'(\eta))$ 
is isotopic to $(D^2 \times S^1, \, \xi_{l'}'(\eta'))$. We have three cases here: when $\infty$ is 
a border between
continued fraction blocks, when $\infty$ is a border between basic slices but not 
between continued fraction blocks, and when $\infty$ is not a border between basic 
slices. 

\medskip
{\bf Case 1}\qua When $\infty$ is a border between continued fraction 
blocks in the sequence $s_0 > \ldots > s_{\frac{n-1}{2n}} > \infty > s_{\frac{m-1}{2m}}' > \ldots > s_0'$, 
it follows from 
the classification theorem for tight contact structures on solid
tori that $(D^2 \times S^1, \, \xi_l'(\eta))$ is isotopic to $(D^2 \times S^1, \, \xi_{l'}'(\eta'))$ if and 
only if $l=l'$, and $\eta$ is isotopic to $\eta'$.

\medskip
{\bf Case 2}\qua The condition for $\infty$ to be a border between basic slices but not
 a border between continued fraction blocks in the minimal basic slices 
decomposition of $(D^2 \times S^1, \, \xi_l'(\eta))$ is that the  slopes $s_{\frac{n-1}{2n}}$ and 
$s_{\frac{m-1}{2m}}'$ are represented by shortest integer vectors $v_{-1}, v_0, v_1$
such that $(v_{-1}, v_0)$ and $(v_0, v_1)$ are integer bases, and $|\det (v_{-1}, 
v_{1})|=2$. This condition is satisfied if and only if $e_0=2$ and $\frac{p}{q} =1$.
In this case the basic slices belonging to the outermost continued fraction 
block of $\eta$, which has boundary slopes $\infty$ and $- \frac{1}{d_0+1}$, and the 
basic slice $(T^2 \times [0, \frac 12 ], \, \zeta_l'')$, with boundary slopes $-1$ and $\infty$  
form a unique continued fraction block in $(D^2 \times S^1, \, \xi_l'(\eta))$, therefore their 
signs can be shuffled. The shuffle can occur in $(M', \, \xi_l'(\eta))$  as well, 
because by Lemma \ref{modello} there is a contact embedding of 
$(D^2 \times S^1, \, \xi_l'(\eta))$ in $(M', \, \xi_l'(\eta))$. 
 We conclude that, when $M'$ has boundary slopes $-1$ and $1$, $(M', \, \xi_l(\eta))$ is isotopic to $(M', \, \xi_{l'}'(\eta'))$ 
 if and only if one of the following holds: either $l=l'$ and $\eta$ is isotopic 
to $\eta'$, or $l'=-l$ and $r_0'=r_0+l$.

\medskip
{\bf Case 3}\qua If $e_0=1$, $s_1= \frac{p}{q}=1$ and $\xi_l'(\eta)$ is tight, then
$(T^2 \times [0, \frac{1}{2}], \, \zeta_l'')$ glues with
the outermost continued fraction block of $\eta$ to give a basic slice with 
boundary slopes $0$ and $- \frac{1}{d_0+1}$. This basic slice forms a continued 
fraction block with the basic slices belonging to the second outermost 
continued fraction block in $\eta$, which has boundary slopes $- \frac{1}{d_0+1}$
and $- \frac{d_1+1}{d_0 d_1 -1}$, therefore their signs can be shuffled. Again,
by Lemma \ref{modello}, the same result holds on $M'$, so we conclude that
$(M', \, \xi_l'(\eta))$ is isotopic to $(M', \, \xi_{l'}'(\eta'))$, when $M'$ has boundary
slopes $0$ and $1$, if and only if one of the 
following holds: either  $l=l'$ and $\eta$ is isotopic 
to $\eta'$, or $l'=-l$, $r_0'=-r_0$, and $r_1'=r_1+l$.
\end{proof}

Now we analyse the states with infinite boundary slopes.
Let $A= S^1 \times [0,1]$ be an annulus and let $\Gamma_A$ be a multicurve on $A$ which 
closes to a homotopically trivial closed curve in $T^2$ if we identify 
$S^1 \times \{ 0 \}$ to $S^1 \times \{ 1 \}$. Let $(M', \, \xi_{\Gamma_A}(\eta))$ be the contact 
manifold  obtained by 
contact $(- \frac{1}{r})$--surgery on a vertical Legendrian curve with twisting 
number $0$ in $(T^2 \times I, \, \xi_{\Gamma_A})$.

\begin{lemma} \label{stato2}
Let $\Gamma_A$ be an abstract dividing set on $A=S^1 \times I$ without homotopically 
trivial closed
curves. Suppose that $\Gamma_A$ closes to an abstract dividing set on $T^2$ consisting
of a unique homotopically trivial closed curve 
if we identify $S^1 \times \{ 0 \}$ and $S^1 \times \{ 1 \}$. If $\Gamma_A'$ is another multicurve on 
$A$ and $\eta'$ a tight contact structure on $D^2 \times S^1$ such that $(M', \, \xi_{\Gamma_A}(\eta))$
is isotopic to $(M', \, \xi_{\Gamma_A'}(\eta'))$, then $\eta$ is isotopic to $\eta'$ and $\Gamma_A$ is 
diffeomorphic to $\Gamma_A'$.
\end{lemma}
\begin{proof}
By Lemma \ref{confronto}, $\eta$ is isotopic to $\eta'$ and $\widehat{\Gamma}_A$ is 
diffeomorphic to $\widehat{\Gamma}_A'$. If $\Gamma_A$ contains an arc with endpoints on 
different boundary components, then $\widehat{\Gamma}_A= \Gamma_A$ and $\widehat{\Gamma}_A'= \Gamma_A'$, 
so we are done. If this is not the case, $\Gamma_A$ consists 
of arcs with both endpoints on the same side, so $\widehat{\Gamma}_A= \Gamma_A$ and $\Gamma_A'$ 
differs from $\Gamma_A$ by a number of closed curves.

We take a multicurve $\Gamma_B$ on $B=S^1 \times [1,2]$ consisting of arcs with both 
endpoints on the same side so that, after identifying $S^1 \times \{ 0 \}$ with 
$S^1 \times \{ 2 \}$, $\Gamma_A \cup \Gamma_B$ closes to some homotopically 
nontrivial 
curves in $T^2$. Then, gluing $(T^2 \times [1,2], \, \xi_{\Gamma_B})$ to $(M', \, \xi_{\Gamma_A}(\eta)) \cong 
(M', \, \xi_{\Gamma_A'}(\eta'))$ yields a generic tight contact structure on $M(e_0,r)$. 
By Theorem \ref{paperino}, the closure of $\Gamma_A \cup \Gamma_B$ is isotopic 
to the closure $\Gamma'_A \cup \Gamma_B$, in particular they have the same 
number of components. This implies that $\Gamma_A$ and $\Gamma'_A$ contain
the same number of closed curves, then  $\Gamma_A$ is diffeomorphic to
$\Gamma_A'$ by Lemma \ref{confronto}.
\end{proof}

\subsection{Analysis of the transitions}
\begin{theorem}\label{eccezionali}
If $e_0 \geq 2$, then $(M(e_0,r), \, \xi_l(\eta))$ is tight for any $\eta$ and any $l$. If 
$e_0=1$, then $(M(e_0,r), \, \xi_l(\eta))$ is tight if and only if $r_0= \frac{ld_0}{2}$.
Moreover, $\xi_{l'}(\eta')$ is isotopic to $\xi_l(\eta)$ if and only if either $l=l'$ and $\eta$
 is isotopic to $\eta'$, or 
\begin{itemize}
\item $l'=-l$ and $r_0'= r_0+l$ when $e_0=2$,
\item  $l'=-l$, $r_0' =-r_0$ and $r_1'= r_1+l$ when $e_0=1$.
\end{itemize}
\end{theorem}

\begin{proof}
Let $W_0 \subset M$ be a convex incompressible vertical torus
with infinite slope and $\# \Gamma_{W_0}=2$ such that the initial state 
$(M \setminus W_0, \, \xi_{l}(\eta)|_{M \setminus W_0}) \cong (M', \, \xi_{l_0}'(\eta_0))$ is contactomorphic  to 
$(M', \, \xi_l'(\eta))$. If $e_0=1$ and $r_0 \neq \frac{ld_0}{2}$, then there is a transition 
from $W_0$ to $W_1$ which brings us to a state $(M \setminus W_1, \, \xi_l(\eta)|_{M \setminus W_1}) \cong 
(M', \, \xi_{l_1}'(\eta_1))$ with boundary slopes $0$ and $1$ which is overtwisted by 
Lemma \ref{stato1}, therefore $(M(e_0,r), \, \xi_l(\eta))$ is overtwisted. In the rest 
of the proof we will suppose either $e_0>1$ or $r_0 = \frac{ld_0}{2}$. 

By induction, we assume we have reached a state $(M \setminus W_i, \, \xi_l(\eta) |_{M \setminus W_i})$ of 
one of the following kinds.
\begin{enumerate}
\item  $(M \setminus W_i, \, \xi_l(\eta)|_{M \setminus W_i})$ is contactomorphic to $(M', \, \xi_{l_i}(\eta_i))$ 
with boundary slopes $\frac{p_i}{q_i} -e_0$ and $\frac{p_i}{q_i}$. If $e_0>2$, then $l_i=l$
and $\eta_i$ is isotopic to $\eta$. If $e_0=2$, then, either $l_i=l$ and $\eta_i$ is isotopic 
to $\eta$, or $l_i=-l$ and $r_0^{i} = r_0+l$. If $e_0=1$, then, either $l_i=l$ and $\eta_i$ is 
isotopic to $\eta$, or $l_i=-l$,  $r_0^{i}=-r_0$ and $r_1^{i}= r_1+l$.
Here $(r_0^i, \ldots ,r_n^i)$ denote the invariants determining $\eta_i$ and $(r_0, \ldots ,r_n)$ 
denote the invariants determining $\eta$.

\item $(M \setminus W_i, \, \xi_l(\eta) |_{M \setminus W_i})$ is contactomorphic 
to the tight contact manifold $(M', \, \xi_{\Gamma_{\Sigma_i}}(\eta))$ obtained by
contact $(- \frac 1r )$--surgery on the $S^1$--invariant contact manifold 
$(T^2 \times I, \, \xi_{\Gamma_{\Sigma_i}})$. Here $\Gamma_{\Sigma_i}$ is a 
multicurve on $S^1 \times I$ which closes to a multicuve  
$\Gamma_{\overline{\Sigma}_i}$ on $T^2$ consisting of a unique homotopically trivial closed 
curve if we identify $S^1 \times \{ 0 \}$ to $S^1 \times \{ 1 \}$.

We denote by $\chi_+(T^2 \setminus \Gamma_{\overline{\Sigma}_i})$ the Euler 
characteristic of the positive region of $T^2 \setminus \Gamma_{\overline{\Sigma}_i}$ and by 
$\chi_-(T^2 \setminus \Gamma_{\overline{\Sigma}_i})$ the Euler 
characteristic of the negative region of $T^2 \setminus \Gamma_{\overline{\Sigma}_i}$.
We have here three cases.
\begin{itemize}
\item If $e_0>2$, then $\chi_+(T^2 \setminus \Gamma_{\overline{\Sigma}_i}) - \chi_-(T^2 \setminus \Gamma_{\overline{\Sigma}_i})=l$ 
and $\eta_i$ is isotopic to $\eta$. 
\item If $e_0=2$, then, either 
$\chi_+(T^2 \setminus \Gamma_{\overline{\Sigma}_i}) - \chi_-(T^2 \setminus \Gamma_{\overline{\Sigma}_i})=l$ 
and $\eta_i$ is isotopic to $\eta$, or 
$\chi_+(T^2 \setminus \Gamma_{\overline{\Sigma}_i}) - \chi_-(T^2 \setminus \Gamma_{\overline{\Sigma}_i})=-l$ and
$r_0^{i} = r_0+l$. 
\item If $e_0=1$, then, either 
$\chi_+(T^2 \setminus \Gamma_{\overline{\Sigma}_i}) - \chi_-(T^2 \setminus \Gamma_{\overline{\Sigma}_i})=l$ and $\eta_i$ 
is isotopic to $\eta$, or 
$\chi_+(T^2 \setminus \Gamma_{\overline{\Sigma}_i}) - \chi_-(T^2 \setminus \Gamma_{\overline{\Sigma}_i})=-l$, 
$r_0^{i}=-r_0$ and $r_1^{i}= r_1+l$.
\end{itemize}
\end{enumerate} 

A transition from the state $M \setminus W_i$
to the state $M \setminus W_{i+1}$ consists of taking a layer $N_i \cong T^2 \times [\frac{1}{2}, 1]
\subset M \setminus W_i$ with boundary $W_{i+1} \cup W_i$, and
moving it from the front to the back, or vice versa. We only consider the case 
when $N_i$ is a front layer. When $N_i$ is a back layer the proof is completely 
analogous. There are two cases, depending on whether the transition changes
the boundary slopes or the division number of the boundary. 

\medskip
{\bf Case 1}\qua This case corresponds to state transitions from 
$(M \setminus W_i, \, \xi_l(\eta)|_{M \setminus W_i})$ to $(M \setminus W_{i+1}, \, \xi_l(\eta)|_{M \setminus W_{i+1}})$ such that $W_i$ 
and $W_{i+1}$ are minimal and at least one of them has finite slope. 
Suppose that $(M \setminus W_i, \, \xi_l(\eta)|_{M \setminus W_i})$ is contactomorphic to
$(M', \xi_{l_i}'(\eta_i))$ and the transition changes the boundary slopes of the state
from $\frac{p_i}{q_i} -e_0$ and $\frac{p_i}{q_i}$ to $\frac{p_{i+1}}{q_{i+1}} -e_0$ and 
$\frac{p_{i+1}}{q_{i+1}}$.
We can assume that the interior of $N_i$ contains no tori with infinite slope.  
If this is not the case, we split the transition in 
two parts, therefore we can assume that $W_{i+1} \subset M \setminus W_i$ has slope 
$\frac{p_{i+1}}{q_{i+1}} \in (\frac{p_i}{q_i}, \infty ]$. We isotope the Seifert fibration
on $M \setminus W_i$ so that $W_{i+1}$ is a fibred torus and the singular fibre $F_i$ is 
contained in $M \setminus N_i$.

By Lemma \ref{infinito} there is a vertical Legendrian curve with twisting 
number $0$ in $M \setminus N_i$. Using such curve we can find a neighbourhood 
$V_{i+1} \subset M \setminus N_i$ of the singular fibre such that $- \partial (M \setminus V_{i+1})$ has infinite 
slope by arguing as in the proof of Lemma \ref{decomposizione}. The 
contact structures $\xi_{l_i}'(\eta_i)|_{V_{i+1}} = \eta_{i+1}$ and $\xi_{l_i}'(\eta_i)|_{M \setminus (W_i \cup V_{i+1})} = 
\xi_{l_{i+1}}'$ are determined by $\eta_i$ and $l_i$ as described in Lemma \ref{stato1}. In 
particular, $\xi_{l_i}'(\eta_i)|_{N_i} = \xi_{l_{i+1}}'|_{N_i}$ is a minimally twisting tight 
contact structure with relative Euler 
class $ \pm \left ( \begin{array}{c} q_i-q_{i+1} \\ p_i-p_{i+1} \end{array} \right )$, with 
the sign depending on $l_{i+1}$. Moving $N_i$ to the back, its relative Euler class 
becomes $\pm \left ( \begin{array}{c} q_i-q_{i+1} \\ p_i-p_{i+1}-e_0(q_i-q_{i+1}) \end{array} 
\right )$, so $\xi_l(\eta)|_{(M \setminus W_{i+1} \cup V)}$ 
is contactomorphic to  $\xi_{l_{i+1}}'(\eta_{i+1})$ with boundary slopes 
$\frac{p_{i+1}}{q_{i+1}} -e_0$ and $\frac{p_{i+1}}{q_{i+1}}$. 
This 
proves that any admissible transition transforms a state of the type described 
in the inductive assumption to another state of the same type.

\medskip
{\bf Case 2}\qua  This case corresponds to state transitions from 
$(M \setminus W_i, \, \xi_l(\eta)|_{M \setminus W_i})$ to $(M \setminus W_{i+1}, \, \xi_l(\eta)|_{M \setminus W_{i+1}})$ such that $W_i$ 
and $W_{i+1}$ have both infinite slope.
Suppose that $(M \setminus W_i, \, \xi_l(\eta)|_{M \setminus W_i})$ is contactomorphic to
$(M', \, \xi_{\Gamma_{\Sigma_i}}(\eta_i))$  and $\# \Gamma_{W_{i+1}}= \# \Gamma_{W_i} \pm 1$. We isotope the Seifert 
fibration $M' \to S^1 \times I$ so that the singular fibre is contained in 
$M' \setminus N_i$ and, fixed a neighbourhood $V_{i+1}$ of the singular fibre
 so that $- \partial (M \setminus V_{i+1})$ has infinite slope, the restrictions of the fibration 
to $N_i$ and $M' \setminus (N_i \cup V_{i+1})$ are $S^1$--bundles. Let $\Sigma$ be a pair of pants and 
let $\overline{\Sigma}$ be a punctured torus obtained by identifying two of the 
boundary components of $\Sigma$. Let $\sigma_i\co \overline{\Sigma} \to M \setminus V_{i+1}$ be a section 
so that $\Sigma_i' = \sigma_i(\Sigma) \subset M'$ is a convex, $\# \Gamma$--minimising surface with 
Legendrian boundary, and $\Sigma_i' \cap W_{i+1}$ is a Legendrian curve. We call 
$\overline{\Sigma}_i' = \sigma_i(\overline{\Sigma}) \subset M \setminus V_{i+1}$.
$(M \setminus W_i, \, \xi_l(\eta)|_{M \setminus W_i})$ is contactomorphic to $(M', \, \xi_{\Gamma_{\Sigma_i'}}(\eta_{i+1}))$, 
therefore, by Lemma \ref{stato2},
$\eta_{i+1}$ is isotopic to $\eta_i$ and $\Gamma_{\overline{\Sigma}_i'}$ is isotopic to 
$\Gamma_{\overline{\Sigma}_i}$. 
If we define $\Sigma_{i+1} \subset (M \setminus (W_{i+1} \cup V_{i+1})$ as $\Sigma_{i+1} = \overline{\Sigma}_i' \setminus W_{i+1}$,  
then $\Gamma_{\overline{\Sigma}_{i+1}}= \Gamma_{\overline{\Sigma}_i'}$ and 
$(M \setminus W_{i+1}, \, \xi_l(\eta)|_{M \setminus W_{i+1}})$ is contactomorphic to 
$(M', \, \xi_{\Gamma_{\Sigma_{i+1}}}(\eta_{i+1}))$.
\end{proof}
As proved in Theorem \ref{exceptionality}, the tight contact structures 
considered in this section become overtwisted after lifting to any finite 
covering of $M(e_0,r)$ induced by a covering of $T^2$, so they are the 
exceptional tight contact structure of Theorem \ref{numerouno}.
The following corollary gives the number of the exceptional tight contact 
structures on $M(e_0,r)$.
\begin{cor} When $e_0 >0$ the number of exceptional tight contact structures 
on $M(e_0,r)$ is finite and positive. It is
\begin{itemize}
\item $2|d_0 (d_1+1) \ldots (d_k+1)|$ if $e_0>2$,
\item $|(d_0-1)(d_1+1) \ldots (d_k+1))|$ if $e_0=2$,
\item $|d_1(d_2+1) \ldots (d_k+1)|$ if $e_0=1$.
\end{itemize}
The last expression has to be interpreted as $2$ when $- \frac{1}{r}=d_0 \in 
\Z$.
\end{cor}
\begin{proof}
By Theorem \ref{eccezionali} and \ref{solidi},  when $e_0 >2$
for any $l$ there are $|d_0 (d_1+1) \ldots (d_k+1)|$ choices for $\eta$ and $2$ choices for
the background, and all choices give distinct tight contact structures, 
therefore the total number of exceptional tight contact structures on $M(e_0,r)$ is 
$2|d_0 (d_1+1) \ldots (d_k+1)|$. 

When $e_0=2$, there are $|d_0 (d_1+1) \ldots (d_k+1)|$ choices for 
$\eta$ for any choice of the background, but not all choices give distinct tight 
contact structures. In fact, any exceptional tight contact structure with
$l=-2$ and $r_0 > d_0$ is isotopic to a tight contact structure with $l=2$, 
therefore we count $|d_0 (d_1+1) \ldots (d_k+1)|$ distinct exceptional tight contact
structures with $l=2$ and $|(d_1+1) \ldots (d_k+1)|$ distinct tight contact structures
with $l=-2$, namely those obtained from $\eta$ with $r_0=d_0$. The total number of 
distinct exceptional tight contact structures up to isotopy on $M(e_0,r)$ with 
$e_0=2$ is therefore $|(d_0-1)(d_1+1) \ldots (d_k+1)|$. 

If $e_0=1$ and $\frac{1}{d} \not \in \Z$, then for any choice of the background
there are $|(d_1+1) \ldots (d_k+1)|$ choices for $\eta$, but not all choices give 
distinct tight contact structures. In fact, we have $|(d_1+1) \ldots (d_k+1)|$ distinct
exceptional tight contact structures with $l=2$ up to isotopy, and 
$|(d_2+1) \ldots (d_k+1)|$ distinct exceptional tight contact structures with $l=-2$
which have not already been counted, namely the ones with $r_i=d_i$.
The total number of exceptional tight contact structures on $M(e_0,r)$ with
$e_0=1$ is therefore $|d_1(d_2+1) \ldots (d_k+1)|$. If $e_0=1$ and $\frac{1}{r} \in \Z$, 
for any $l$ there is only one possibility for $\eta$, and different choices for the
background produce non isotopic tight contact structures, therefore the total 
of exceptional tight contact structures on $M(e_0,r)$ with $e_0=1$ and 
$\frac 1r \in \Z$ is $2$.
\end{proof}

The exceptional tight contact structures on $M(e_0,r)$ are negative contact 
surgeries on the exceptional tight contact structures on $T(e_0)$ when $e_0>1$.
On the contrary, there are no exceptional tight contact structures on $T(1)$,
 and the exceptional tight contact structures on $M(1 ,r)$ are negative
contact surgeries on an overtwisted contact structure on $T(1)$. On $T(2)$ 
there is only one exceptional tight contact structure up to isotopy, therefore
 the two backgrounds
extend to isotopic tight contact structure on $T(2)$. This reflects the fact
that $T(2)$ with the exceptional tight contact structure contains two non 
Legendrian isotopic 
 vertical Legendrian curves with twisting number $0$, and negative contact 
surgeries
with the same surgery data on such curves yield different contact 
manifolds. The shuffling between the background and the surgery data when
$e_0=2$ shows that suitably stabilisations of the two non isotopic vertical 
Legendrian curves with twisting number $0$ become Legendrian isotopic.
 
%%%%%%%%%%%%%%%%%%%%%%%%%%%%%%%%%%%%%%%%%%%%%%%%%%%%%%%%%%%%%%%%%%%%%%%%%%%%%%%

\Addresses\recd
\end{document}